\documentclass[preprint]{elsarticle}
\journal{Journal of Computational Physics}

\usepackage{mathrsfs}
\usepackage{amsmath}
\usepackage{amsfonts}
\usepackage{amssymb}
\usepackage[utf8]{inputenc}
\usepackage[T1]{fontenc}

\usepackage{graphicx}
\usepackage{xcolor, colortbl}
\usepackage{array,color}
\usepackage{arydshln}
\usepackage{threeparttable}
\usepackage{multirow}
\usepackage{float}
\usepackage{bm}

\usepackage[ruled,vlined]{algorithm2e}
\usepackage{algorithmic}
\usepackage{hhline}

\usepackage{upgreek}
\usepackage{tikz}
\usetikzlibrary{shadows}
\usetikzlibrary{shapes}
\usepackage{soulutf8}

\newtheorem{thm}{Theorem}

\newtheorem{proposition}[thm]{Proposition}
\newdefinition{rmk}{Remark}
\newproof{pf}{Proof}

\begin{document}

\begin{frontmatter}

\title{Schwarz waveform relaxation method for one dimensional Schr{\"o}dinger equation with general potential}

\author[adresseCB]{C. Besse}
\ead{christophe.besse@math.univ-toulouse.fr}

\author[adresseXF]{F. Xing}
\ead{feng.xing@unice.fr}

\address[adresseCB]{Institut de Math{\'e}matiques de Toulouse UMR5219,
        Universit\'e de Toulouse; CNRS,
        UPS IMT, F-31062 Toulouse Cedex 9, France.}
\address[adresseXF]{Maison de la Simulation, CEA Saclay France \& Laboratoire Paul Painlev{\'e}, Universit{\'e} Lille Nord de France.}

\begin{abstract}
  In this paper, we apply the Schwarz Waveform Relaxation (SWR) method to the one dimensional Schr{\"o}dinger equation with a general linear or a nonlinear potential. We propose a new algorithm for the Schr{\"o}dinger equation with time independent linear potential, which is robust and scalable up to 500 subdomains. It reduces significantly computation time compared with the classical algorithms. Concerning the case of time dependent linear potential or the nonlinear potential, we use a preprocessed linear operator for the zero potential case as preconditioner which lead to a preconditioned algorithm. This ensures high scalability. Besides, some newly constructed absorbing boundary conditions are used as the transmission condition and compared numerically. 
\end{abstract}



\end{frontmatter}

\section{Introduction}

Schwarz waveform relaxation method (SWR) is one class of the domain decomposition methods for time dependent partial differential equations. The time-space domain is decomposed into subdomains. The solution is computed on each subdomain for whole time interval and exchange the time-space boundary value. Some articles are devoted to this method for linear Schr{\"o}dinger equation \cite{Halpern2010_sch, Halpern2006_sch}, advection reaction diffusion equations \cite{Caetano2010, Gander2007_diffusion, Hoang2013}, wave equations \cite{Gander2003_wave, Halpern2009_wave} and Maxwell's equation \cite{Dolean2009max}.

This paper deals with the SWR method without overlap for the one dimensional Schr{\"o}dinger equation defined on a bounded spatial domain $(a_0,b_0)$, $a_0, b_0 \in \mathbb{R}$ and $t \in (0,T)$. The Schr{\"o}dinger equation with homogeneous Neumann boundary condition reads
\begin{equation}
  \label{Schequ}
  \left\{\begin{array}{ll}
      \mathscr{L}u := (i\partial_t + \partial_{xx} + \mathscr{V} )u = 0, \ (t,x)\in (0,T)\times (a_0,b_0),\\
      u(0,x) = u_0(x), \ x \in (a_0,b_0),\\
      \partial_\mathbf{n} u (t,x) = 0, \ x=a_0,b_0,
    \end{array} \right.
\end{equation}
where $\mathscr{L}$ is the Schr{\"o}dinger operator, $\partial_\mathbf{n}$ is the normal directive, the initial value $u_0 \in L^2(\mathbb{R})$ and $\mathscr{V}$ is a real potential. We consider both linear and nonlinear potentials: 
\begin{enumerate}
\item $\mathscr{V}=V(t,x)$,
\item $\mathscr{V}=f(u)$, ex. $\mathscr{V}=|u|^2$.
\end{enumerate}
%
In order to perform domain decomposition method, the time-space domain $(0,T) \times (a_0,b_0)$ is decomposed into $N$ subdomains $\Uptheta_j = (0,T)\times \Omega_j$, $\Omega_j=(a_j,b_j)$ without overlap as shown in Figure \ref{chp1_sub3} for $N=3$. 
\definecolor{col1}{rgb}{0.8,1,0.8}
\definecolor{col2}{rgb}{1,0.8,0.8}
\definecolor{col3}{rgb}{0.8,0.8,1}
\begin{figure}[H]
  \centering
  \begin{tikzpicture}
    \draw[>=stealth,->] (-0.5,0) -- (8,0);
    \draw (8,0) node[right] {$x$};
    \draw [>=stealth,->] (0,0) -- (0,2.5);
    \draw (0,2.5) node[above] {$t$};
    \draw [color=gray,fill=col1] (0.2,0) rectangle (2,2);
    \draw [color=gray,fill=col2] (2,0) rectangle (4.6,2);
    \draw [color=gray,fill=col3] (4.6,0) rectangle (7.5,2);
    
    \draw (0.2,-0.1) node[below] {$a_0=a_1$};
    \draw (2,0) node[below] {$b_1=a_2$};
    \draw (4.6,0) node[below] {$b_2=a_3$};
    \draw (7.5,0) node[below] {$b_2=b_0$};
    \draw (0,2) node[left] {$T$};
    \draw (1.1,1) node[scale=0.7] {$(0,T)\times \Omega_1$};
    \draw (3.3,1) node[scale=0.7] {$(0,T)\times \Omega_2$};
    \draw (6.05,1) node[scale=0.7] {$(0,T)\times \Omega_3$};

    \draw[>=stealth,->] (2,1.8) -- (1.7,1.8);
    \draw (1.7,1.8) node[below,scale=0.7] {$\mathbf{n}_2$};

    \draw[>=stealth,->] (4.6,1.8) -- (4.9,1.8);
    \draw (4.9,1.8) node[below,scale=0.7] {$\mathbf{n}_2$};
  \end{tikzpicture}
  \caption{Domain decomposition without overlap, $N=3$.}
  \label{chp1_sub3}
\end{figure}
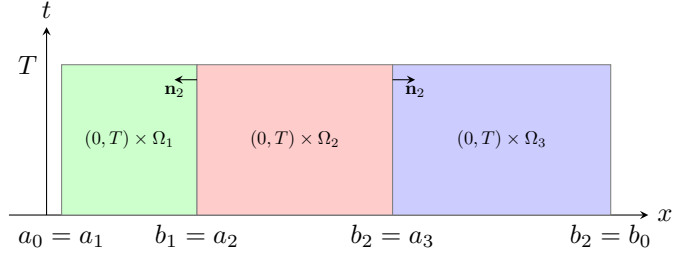
The classical SWR algorithm consists in applying the sequence of iterations for $j=2,3,...,N-1$
\begin{equation} 
  \label{Algo}
  \left\{ \begin{array}{ll}
      \mathscr{L} u_j^{k+1} = 0, &(t,x) \in \Uptheta_j,\\[1mm]
      u^{k+1}_j(0,x) = u_0(x), &x \in \Omega_j, \\[1mm]
      B_j u_j^{k+1} = B_j u_{j-1}^{k}, &x = a_j, \\[1mm]
      B_j u_j^{k+1} = B_j u_{j+1}^{k}, &x = b_j.
    \end{array} \right.
\end{equation}
The two extremal subdomains require special treatment since the Neumann boundary condition is imposed in \eqref{Schequ} at the points $a_0$ and $b_0$.
\begin{align*}
  \left\{ \begin{array}{ll}
      \mathscr{L} u_1^{k+1} = 0, (t,x) \in \Uptheta_1,\\[1mm]
      u^{k+1}_1(0,x) = u_0(x), x \in \Omega_1, \\[1mm]
      \partial_{\mathbf{n}_1} u_1^{k+1} = 0, x = a_1, \\[1mm]
      B_1 u_1^{k+1} = B_1 u_2^{k} , x = b_1, 
    \end{array} \right.&&
  \left\{ \begin{array}{ll}
      \mathscr{L} u_N^{k+1} = 0, (t,x) \in \Uptheta_N,\\[1mm]
      u^{k+1}_N(0,x) = u_0(x), x \in \Omega_N, \\[1mm]
      B_N u_N^{k+1} = B_N u_{N-1}^{k}, x = a_N, \\[1mm]
      \partial_{\mathbf{n}_N} u_N^{k+1} = 0, x = b_N.
    \end{array} \right.
\end{align*}
The notation $u_j^{k}$ denotes the solution on subdomain $\Uptheta_j=(0,T)\times (a_j,b_j)$ at iteration $k=0,1,2,...$ of the SWR algorithm. The boundary information is transmitted with adjacent subdomains $\Uptheta_{j-1}$ and $\Uptheta_{j+1}$ through the transmission operators $B_j$.

The transmission condition is one of the key issues for this method. For the linear Schr{\"o}dinger equation, the SWR method with or without overlap is introduced and analyzed by Halpern and Szeftel in \cite{Halpern2010_sch}. For the decomposition without overlap, if $\mathscr{V}$ is a constant, they use an optimal transmission condition given by the underlying transparent boundary condition. However, the transparent boundary condition is not always available for a variable potential. Robin transmission condition and quasi-optimal transmission condition are therefore used and are named as optimized Schwarz waveform relaxation algorithm and quasi-optimal Schwarz waveform relaxation algorithm respectively. In both cases, the transmission operator is written as
\begin{equation} \label{S1S2}
  B_j = \partial_{\mathbf{n}_j} + S_j,
\end{equation}
where the operator $S_j$ is
\begin{displaymath}
    \text{Robin}: \ S_j = - ip, \ p \in \mathbb{R}^{+}, \quad 
    \text{Quasi-optimal}:  \ S_j  = \sqrt{-i\partial_t - V|_{a_j,b_j}},
\end{displaymath}
and $\mathbf{n}_j$ denotes the outwardly unit normal vector at $a_j$ or $b_j$.
Recently, Antoine, Lorin and Bandrauk \cite{Antoine2014} consider the general Schr{\"o}dinger equation. On the interface between subdomains, they propose to use recent absorbing conditions as transmission condition, which is also an idea that we follow in this paper.

In recent years, some absorbing operators for one dimensional Schr{\"o}dinger equation have been constructed by using some adaptations of pseudo-differential techniques \cite{Antoine2006,Antoine2009a,Antoine2011,Antoine2009}. We use them here as the transmission operators in \eqref{S1S2} and expect to get good convergence properties. 

We are also interested in this article about the effectivness of the method on parallel computers. Another import issue for the method is therefore the scalability. As we know, without additional considerations, the more subdomains are used to decomposed $(a_0,b_0)$, the more iterations are required for SWR algorithm to reach convergence. Thus, the total computation time could hardly decrease significantly. In this paper, we propose two solutions: a new scalable algorithm if the potential is independent of time and a preconditioned algorithm for general potentials.

This paper is organized as follows. In section 2, we present the transmission conditions which are used in this paper for the classical SWR algorithm, and the discretization that plays an important role for the analyses of the interface problem in Section 3. In Section 4 and 5, we present the new algorithm for time independent linear potential and the preconditioned algorithm for general potentials. Some numerical results are shown in Section 6. Finally, we draw a conclusion in the last section.

\section{SWR algorithm and discretization}
\label{Sec_SWR}

\subsection{Transmission conditions}

The transmission conditions on boundary points $a_j$ and $b_j$ are given thanks to the relation
\begin{equation} 
\label{chp1_S1S2}
  B_j = \partial_{\mathbf{n}_j} + S_j,
\end{equation}
where the operators $S_j$ could take different forms. Besides the Robin transmission condition, we propose in this paper to use the operators $S_j$ coming from the artificial boundary conditions for \eqref{Schequ} defined in \cite{Antoine2009, Antoine2009a, Antoine2011, Klein2010these} for a linear or nonlinear potential $\mathscr{V}(t,x,u)$. The authors propose three families of conditions written as
\begin{displaymath}
\partial_{\mathbf{n}} u+ S_l^M u=0,
\end{displaymath}
on the boundary of considered computation domain, $M$ denotes the order of the artificial boundary conditions. We index by $l$ these families of boundary conditions: $l=0$ for potential strategy, $l=1$ for gauge change strategy and $l=2$ for Pad{\'e} approximation strategy. We recall here the definition of operators $S_l^M$ for the different strategies.

\vbox{}
\noindent \textbf{Potential strategy $l=0$} (\cite{Antoine2009})
\begin{align*}
    \mathrm{Order} \ 2: & \quad S_0^2 = e^{-i\frac{\pi}{4}}\partial_t^{1/2},\\
    \mathrm{Order} \ 3: & \quad S_0^3 = S_0^2 - e^{i\frac{\pi}{4}} \frac{\mathscr{V}}{2} I_t^{1/2}, \\
    \mathrm{Order} \ 4: &\quad  S_0^4 = S_0^3 - i\frac{\partial_\mathbf{n} \mathscr{V}}{4} I_t,
\end{align*}
where the fractional half-order derivative operator $\partial_t^{1/2}$ applied to a function $h$ is defined by
\begin{displaymath}
  \partial_t^{1/2} h(t) = \frac{1}{\sqrt{\pi}} \partial_t \int_0^t \frac{h(s)}{\sqrt{t-s}}ds,
\end{displaymath} 
the half-order integration operator $I_t^{1/2}$ and the integration operator are given by
 \begin{displaymath}
   I_t^{1/2} h(t) = \frac{1}{\sqrt{\pi}} \int_0^t \frac{h(s)}{\sqrt{t-s}}ds, \ \ I_t h(t) = \int_0^t h(s)ds.
 \end{displaymath}

\vbox{}
\noindent \textbf{Gauge change strategy $l=1$} (\cite{Antoine2009a, Antoine2011})
\begin{align*}
    \mathrm{Order} \ 2: & \quad  S_1^2 = e^{-i\frac{\pi}{4}} e^{i\mathcal{V}(t,x)}\partial_t^{1/2}(e^{-i\mathcal{V}(t,x)} \cdot), \\
    \mathrm{Order} \ 4: & \quad  S_1^4 = S_1^2 - i\mathrm{sgn}(\partial_\mathbf{n} \mathscr{V})\frac{\sqrt{|\partial_\mathbf{n} \mathscr{V}|}}{2}e^{i\mathcal{V}(t,x)} I_t(\frac{\sqrt{|\partial_\mathbf{n} \mathscr{V}|}}{2}e^{-i\mathcal{V}(t,x)} \cdot),
\end{align*}
where $sgn(\cdot)$ is the sign function and 
\begin{displaymath}
  \mathcal{V}(t,x) = \int_0^t  \mathscr{V}(s,x,u(s,x))ds.
\end{displaymath}

\vbox{}
\noindent \textbf{Pad{\'e} approximation strategy $l=2$} (\cite{Antoine2009a, Antoine2011})
\begin{align*}
    \mathrm{Order} \ 2: & \quad  S_2^2 =  -i \sqrt{i\partial_t + \mathscr{V}}, \\ 
    \mathrm{Order} \ 4: & \quad  S_2^4 = S_2^2 + \mathrm{sgn}(\partial_{\mathbf{n}} \mathscr{V}) \frac{\sqrt{|\partial_{\mathbf{n}} \mathscr{V}|}}{2} \big( i\partial_t+ \mathscr{V} \big)^{-1} \Big( \frac{\sqrt{|\partial_{\mathbf{n}} \mathscr{V}|}}{2} \cdot \Big).
\end{align*}

\subsection{Discretization}

The aim of this subsection is to present the discretization of the Schr{\"o}dinger equation with a linear potential $\mathscr{V}=V(t,x)$ or a nonlinear potential $\mathscr{V}=f(u)$.

\subsubsection{Case of linear potential}

First, we describe the discretization of the linear Schr{\"o}dinger equation. We discretize the time interval $(0,T)$ uniformly with $N_T$ intervals and define $\Delta t = T/N_T$ to be  the time step. A semi-discrete approximation adapted to the Schr{\"o}dinger equation on $(0,T) \times (a_j,b_j),j=1,2,...,N$ is given by the semi-discrete Crank-Nicolson scheme
\begin{displaymath}
  i\frac{u_{j,n}^{k} - u_{j,n-1}^{k}}{\Delta t} + \partial_{xx} \frac{u_{j,n}^{k} + u_{j,n-1}^{k}}{2} + \frac{V_{n}+V_{n-1}}{2} \frac{u_{j,n}^{k} + u_{j,n-1}^{k}}{2} = 0, \ 1 \leqslant n \leqslant N_T,
\end{displaymath} 
and $u_{j,0}^k=u(0,x)$ for $x\in (a_j,b_j)$. The unknown function $u_{j,n}^{k}(x)$ is an approximation of the solution $u_j^k(n\Delta t,x)$ to the Schr\"odinger equation at time $t_n=n\Delta t$ on subdomain $\Omega_j$ and at iteration $k$. We define the approximation of the potential $V_n(x)=V(t_n,x)$.

For implementation issue, it is useful to introduce new variables $v_{j,n}^{k} = (u_{j,n}^{k} + u_{j,n-1}^{k})/2$ with $v^{k}_{j,0} = u^{k}_{j,0}$. The scheme could be written as
\begin{equation}
\label{CNS_L}
  2 i\frac{v_{j,n}^{k}}{\Delta t} + \partial_{xx} v_{j,n}^{k} +  W_{n} v_{j,n}^{k} = 2i\frac{u_{j,n-1}^{k}}{\Delta t},
\end{equation}
with $W_n=(V_n+V_{n-1})/2$.
The spatial approximation is realized thanks to a classical $P_1$ finite element method.
%
%
%
The use of transmission condition gives the following boundary conditions for each subdomain 
  \begin{equation}
    \label{condlimit}
    \hspace{0.1cm} \ \left \{ \begin{array}{ll}
        \partial_{\mathbf{n}_j} v^{k}_{j,n} + \overline{S} v^{k}_{j,n} = \partial_{\mathbf{n}_j} v^{k-1}_{j-1,n} + \overline{S} v^{k-1}_{j-1,n},\ x=a_j, \\[2mm]
        \partial_{\mathbf{n}_j} v^{k}_{j,n} + \overline{S} v^{k}_{j,n} = \partial_{\mathbf{n}_j} v^{k-1}_{j+1,n} + \overline{S} v^{k-1}_{j+1,n}, \ x=b_j,
      \end{array} \right. 
  \end{equation}
with special treatments for the two extreme subdomains
\begin{displaymath}
\partial_{\mathbf{n}_1} v^{k}_{1,n} = 0, \ x=a_1, \quad \partial_{\mathbf{n}_N} v^{k}_{N,n} = 0, \ x=b_N,
\end{displaymath}
%
%
where $\overline{S}$ is a semi-discretization of $S$. For each strategy, $\overline{S}$ is given by

\vbox{}
\noindent \noindent \textbf{Potential strategy} $l=0$
%
\begin{align*}
  \mathrm{Order} \quad 2: \quad  \overline{S}_{0}^{2} v_{j,n}^{k} &= e^{-i\pi/4}\sqrt{\frac{2}{\Delta t}}\sum_{s=0}^{n}\beta_{n-s}v^{k}_{j,s}, \\
  \mathrm{Order} \quad 3: \quad  \overline{S}_{0}^{3} v_{j,n}^{k} &= \overline{S}_{0}^{2} v_{j,n}^{k} - e^{i\pi/4} \sqrt{\frac{\Delta t}{2}}\frac{W_{n}}{2}\sum_{s=0}^{n}\alpha_{n-s}v^{k}_{j,s},\\
  \mathrm{Order} \quad 4: \quad  \overline{S}_{0}^{4} v_{j,n}^{k} &= \overline{S}_{0}^{3} v_{j,n}^{k} - i \frac{\partial_{\mathbf{n}_j} W_{n}}{4} \frac{\Delta t}{2}\sum_{s=0}^{n}\gamma_{n-s}v^{k}_{j,s},
\end{align*}
where
\begin{displaymath}
  \begin{split}
    &(\alpha_0,\alpha_1,\alpha_2,\alpha_3,\alpha_4,\alpha_5,...) = (1,1,\frac{1}{2},\frac{1}{2},\frac{3}{8},\frac{3}{8},\frac{3 \cdot 5}{2 \cdot 4 \cdot 6},...), \
    \beta_s = (-1)^s \alpha_s, \forall s \geqslant 0, \\
    &(\gamma_0,\gamma_1,\gamma_2,\gamma_3,...) = (1,2,2,2,...).
  \end{split}
\end{displaymath}

\vbox{}
\noindent \noindent \textbf{Gauge change strategy} $l=1$
\begin{align*}
  \mathrm{Order} \quad 2: \quad  & \overline{S}_{1}^{2} v_{j,n}^{k} = e^{-i\pi/4}e^{i\mathcal{W}_{n}}\sqrt{\frac{2}{\Delta t}}\sum_{s=0}^{n}\beta_{n-s}e^{-i\mathcal{W}_{s}}v_{j,s}^k, \\
  \mathrm{Order} \quad 4: \quad  & \overline{S}_{1}^{4} v_{j,n}^{k} = \overline{S}_{1}^{2} v_{j,n}^{k}\\
  & - i sgn(\partial_{\mathbf{n}_j}W_{n})\frac{\sqrt{|\partial_{\mathbf{n}_j} W_{n}|}}{2}e^{i\mathcal{W}_{n}}\frac{\Delta t}{2} \sum_{s=0}^{n} \gamma_{n-s}\frac{\sqrt{|\partial_{\mathbf{n}_j} W_{s}|}}{2}e^{-i\mathcal{W}_{s}}v_{j,s}^k,
\end{align*}
where $\mathcal{W}_{n} = \frac{\mathcal{V}_{n}+\mathcal{V}_{n-1}}{2}$ and $\mathcal{V}_n(x)=\int_0^{t_n} V(s,x) ds$.

\vbox{}
\noindent \noindent \textbf{Pad{\'e} approximation strategy} $l=2$
\begin{align*}
  \overline{S}_{2}^{2} v_{j,n}^{k} = & -i \big( \sum_{s=0}^{m} a_s^{m} \big) v_{j,n}^{k} + i\sum_{s=1}^{m}a_s^{m}d_s^{m} \frac{1}{\frac{2i}{\Delta t}+W_{n}+d_s^m} v_{j,n}^k\\
  & + i\sum_{s=1}^{m}a_s^{m}d_s^{m} \frac{\frac{2i}{\Delta t}}{\frac{2i}{\Delta t}+W_{n}+d_k^m} \varphi_{j,n-1}^{s}, \\
  \overline{S}_{2}^{4} v_{j,n}^{k} = & \overline{S}_{2}^{2} v_{j,n}^{k} + \frac{\partial_{\mathbf{n}_j} W_{n}}{4}  \frac{1}{\frac{2i}{\Delta t} + W_{n} }  v_{j,n}^{k} \\
  & + sgn(\partial_{\mathbf{n}_j}W_{n})\frac{\sqrt{|\partial_{\mathbf{n}_j}W_{n}|}}{2} \frac{\frac{2i}{\Delta t}}{\frac{2i}{\Delta t} + W_{n} } \psi_{j,n-1},
\end{align*}
where $\varphi_{j,n}^{s},\phi_{j,n},s=1,2,...,m$ are introduced as auxiliary functions
\begin{displaymath}
  \left\{
    \begin{array}{l}
      \varphi_{j,n-\frac{1}{2}}^s = \frac{1}{\frac{2i}{\Delta t} + W_{n}+d_s^m} v_{j,n}^k + \frac{\frac{2i}{\Delta t}}{\frac{2i}{\Delta t} + W_{n}+d_s^m} \varphi_{j,n-1}^s, \ s=1,2,...,m\\
      \varphi_{j,n}^s = 2 \varphi_{j,n-\frac{1}{2}}^s - \varphi_{j,n-1}^s, \\ 
      \varphi_{j,0}^s = 0,
    \end{array}
  \right.
\end{displaymath}
and
\begin{displaymath}
  \left\{
    \begin{array}{l}
      \psi_{n-\frac{1}{2}}  = \frac{\sqrt{|\partial_{\mathbf{n}_j}W_{n}|}}{2} \frac{1}{\frac{2i}{\Delta t} + W_{n} } v_{j,n}^{k} + \frac{ \frac{2i}{\Delta t} }{\frac{2i}{\Delta t} + W_{n} }  \psi_{j,n-1},\\
      \psi_{j,n} = 2 \psi_{j,n-\frac{1}{2}} - \psi_{j,n-1},\\
      \psi_{j,0} = 0.
    \end{array}
  \right.
\end{displaymath}
We also recall here the Robin transmission condition and its approximation
\begin{displaymath}
  S = S_p = -ip,  \quad \overline{S} v_{j,n}^{k} = \overline{S}_p v_{j,n}^{k} = -ip \cdot v_{j,n}^{k}, \quad p \in \mathbb{R}^{+}.
\end{displaymath}

We propose below to rewrite \eqref{condlimit} by using fluxes, which are defined at interfaces by 
  \begin{align*}
    l_{j,n}^{k} = \partial_{\mathbf{n}_j} v^{k}_{j,n}(a_j) + \overline{S} v^{k}_{j,n}(a_j),\
    r_{j,n}^{k} = \partial_{\mathbf{n}_j} v^{k}_{j,n}(b_j) + \overline{S} v^{k}_{j,n}(b_j), \ j=1,2,...,N,
  \end{align*}
with the exception for $l_{1,n}^{k} = r_{N,n}^{k} = 0$. It is obvious that, on each subdomain, the boundary conditions are
  \begin{equation}
     \left \{ \begin{array}{ll}
     \label{TransmissionCondlocal}
        \partial_{\mathbf{n}_j} v^{k}_{j,n} + \overline{S} v^{k}_{j,n} = l_{j,n}^k,\ x=a_j, \\[2mm]
        \partial_{\mathbf{n}_j} v^{k}_{j,n} + \overline{S} v^{k}_{j,n} = r_{j,n}^k, \ x=b_j.
      \end{array} \right. 
  \end{equation}
 For the transmission condition $S_0^2$, $S_0^3$, $S_1^2$ and $S_2^2$ which do not contain 
the normal derivative of potential $W_n$, using \eqref{condlimit}, we have
\begin{align*}
  r_{1,n}^{k} & = \partial_{\mathbf{n}_1} v^{k}_{1,n}(b_1) + \overline{S} v^{k}_{j,n}(b_1) = \partial_{\mathbf{n}_1} v^{k-1}_{2,n}(a_{2}) + \overline{S} v^{k-1}_{2,n}(a_{2}) \\
  & -\Big( \partial_{\mathbf{n}_{2}} v^{k-1}_{2,n}(a_{2}) + \overline{S} v^{k-1}_{2,n}(a_{2}) \Big) + 2\overline{S} v^{k-1}_{2,n}(a_{2}) = -l_{2,n}^{k-1} + 2\overline{S} v^{k-1}_{2,n}(a_{2}).
\end{align*}
The transmission conditions could therefore be rewritten as
  \begin{equation}
   \label{TransmisionCond}
   \left \{
  \begin{array}{ll}
        l_{1,n}^{k} = 0, \ \ l_{j,n}^{k} = -r_{j-1,n}^{k-1} + 2\overline{S} v_{j-1,n}^{k-1}(b_{j-1}), \ j=2,...,N, \\[2mm]
        r_{1,N}^{k} = 0, \  r_{j,n}^{k} = -l_{j+1,n}^{k-1} + 2\overline{S} v_{j+1,n}^{k-1}(a_{j+1}), \ j=1,2,...,N-1.
  \end{array} \right.
  \end{equation} 
Dealing with the transmission conditions $S_0^4$, $S_1^4$ and $S_2^4$, we could also obtain similar formulas to \eqref{TransmisionCond}. We can therefore replace the boundary conditions \eqref{condlimit} for the $N$ local problems \eqref{CNS_L} by \eqref{TransmissionCondlocal} and  fluxes definition \eqref{TransmisionCond}.

Let us denote by $\mathbf{v}_{j,n}^k$ (resp. $\mathbf{u}_{j,n}^k$) the nodal $P_1$ interpolation vector of $v_{j,n}^k$ (resp. $u_{j,n}^k$) with $N_j$ nodes, $\mathbb{M}_j$ the mass matrix, $\mathbb{S}_j$ the stiffness matrix and $\mathbb{M}_{j,W_n}$ the generalized mass matrix with respect to $\int_{a_j}^{b_j} W_n v \phi dx$, $j=1,2,...,N$. Thus, the matrix formulation of the $N$ local problems is given by
\begin{equation}
  \label{AvML}
  (\mathbb{A}_{j,n} - \mathbb{B}_{j,n}) \mathbf{v}_{j,n}^k = \frac{2i}{\Delta t} \mathbb{M}_j\mathbf{u}_{j,n-1}^k + \mathbf{b}_{j,n}^k - Q_j^T (l_{n}^k, r_{n}^k)^T,  
\end{equation}
where $\mathbb{A}_{j,n} = \frac{2i}{\Delta t}\mathbb{M}_j - \mathbb{S}_j + \mathbb{M}_{j,W_{n}}$ and "$\cdot^T$" is the standard notation of the transpose of a matrix or a vector. The restriction matrix $Q_j$ is defined by
\begin{displaymath}
  Q_{j} = \begin{pmatrix}
    1 & 0 & 0 & \cdots & 0 & 0 \\
    0 & 0 & 0 & \cdots & 0 & 1
  \end{pmatrix}
  \in \mathbb{C}^{2 \times N_j}.
\end{displaymath}
$\mathbb{B}_{j,n} \in \mathbb{C}^{N_j\times N_j}$ (resp. $\mathbf{b}_{j,n}^k \in \mathbb{C}^{N_j}$) represent the boundary matrix (resp. vector) associated with the boundary condition at time step $n$, which depends on the transmission condition. The discrete form of the transmission condition \eqref{TransmisionCond} is given by
\begin{equation}
  \label{TransmissionCond_disc}
  \left \{ \begin{array}{ll}
    l_{j,n}^{k} = -r_{j-1,n}^{k} + 2 \widetilde{S} (Q_{j-1,r} \mathbf{v}_{j-1,n}^{k}), \ j=1,2,...,N-1,\\[2mm]
    r_{j,n}^{k} = -l_{j+1,n}^{k} + 2 \widetilde{S} (Q_{j+1,l} \mathbf{v}_{j+1,n}^{k}),\ j=2,3,...,N.
  \end{array} \right.
\end{equation} 
where $Q_{j,l} = (1, 0, \cdots, 0, 0 ) \in \mathbb{C}^{N_j}$, $Q_{j,r} = (0, 0, \cdots, 0, 1) \in \mathbb{C}^{N_j}$. $\widetilde{S}$ is the fully discrete version of $\overline{S}$. For example the transmission condition $S_0^2$ leads to
\begin{align*}
      & \widetilde{S}_{0}^{2} (Q_{j,l} \mathbf{v}_{j,n}^{k}) = e^{-i\pi/4}\sqrt{\frac{2}{\Delta t}} \sum_{s=0}^{n}\beta_{n-s} (Q_{j,l} \mathbf{v}^{k}_{j,s}),\\
      & \widetilde{S}_{0}^{2} (Q_{j,r} \mathbf{v}_{j,n}^{k}) = e^{-i\pi/4}\sqrt{\frac{2}{\Delta t}} \sum_{s=0}^{n}\beta_{n-s} (Q_{j,r} \mathbf{v}^{k}_{j,s}).
\end{align*}

\subsubsection{Case of nonlinear potential}

If the potential is nonlinear $\mathscr{V}=f(u)$, we propose to use the usual scheme developed by Durán- Sanz Serna \cite{DUR2000}
\begin{displaymath}
  i\frac{u_{j,n}^{k} - u_{j,n-1}^{k}}{\Delta t} + \partial_{xx} \frac{u_{j,n}^{k} + u_{j,n-1}^{k}}{2} + f(\frac{u_{j,n}^{k} + u_{j,n-1}^{k}}{2}) \frac{u_{j,n}^{k} + u_{j,n-1}^{k}}{2} = 0,\ 1 \leqslant n \leqslant N_T,
\end{displaymath}
By using the notations defined in the previous subsection, this schema reads as
\begin{equation}
\label{CNS_NL}
  2 i\frac{v_{j,n}^{k}}{\Delta t} + \partial_{xx} v_{j,n}^{k} +
  f(v_{j,n}^{k}) v_{j,n}^{k} = 2i\frac{u_{j,n-1}^{k}}{\Delta t}.
\end{equation}
As in the previous subsection, we use a $P_1$ finite element method to deal with the space variable approximation. 
Since the problem is nonlinear, the computation of $v_{j,n}^k$ is made by a fixed point procedure. At a given time $t=t_n$, we take $\bm{\zeta}_j^{0} = \mathbf{v}^k_{j,n-1}$ and compute the solution $\mathbf{v}^k_{j,n}$ as the limit of the iterative scheme with respect to $s$:
\begin{equation}
  \label{AvMNL}
  \left(\frac{2i}{\Delta t}\mathbb{M}_{j} - \mathbb{S}_{j} - \mathbb{B}_{j,n}\right)
  \bm{\zeta}^{s+1}_{j} \\
  = \frac{2i}{\Delta t} \mathbb{M}_j \mathbf{u}_{j,n-1}^{k} -
  \mathbf{b}_{j,f(\zeta^{s}_j)} + \mathbf{b}_{j,n}^k - Q_j^T (l^{k}_{n},r^{k}_{n})^T,
\end{equation}
where $\mathbf{b}_{j,f(v)}$ is the vector associated with $\int_{a_j}^{b_j} f(v) v \phi dx$. The matrix $\mathbb{B}_{j,n}$ and the vector $\mathbf{b}_{j,n}^k$ depend on the transmission operator. The discrete form of the transmission conditions is similar to \eqref{TransmissionCond_disc} obtained for linear potential.

\section{Interface problem}
\label{Sec_InterfacePb}

The $N$ problems \eqref{AvML} and \eqref{AvMNL} on each subdomain could be written globally. Let us define the global interface vector $g^k$ at iteration $k$ by
\begin{equation*}
g^k = \big( \underbrace{r_{1,1}^k,r_{1,2}^k,...,r_{1,N_T}^k}_{j=1}, \cdots, \underbrace{l_{j,1}^k,...,l_{j,N_T}^k,r_{j,1}^k,...,r_{j,N_T}^k}_j, \cdots,
\underbrace{l_{N,1}^k,l_{N,2}^k,...,l_{N,N_T}^k}_{j=N} \big)^T.
\end{equation*}
Considering the transmission conditions with flux \eqref{TransmissionCond_disc}, it is not hard to see that there exist an operator $\mathcal{R}$ such that
 \begin{equation}
 \label{interfacepbR}
 g^{k+1} = \mathcal{R} g^{k}.
 \end{equation}
The interface operator $\mathcal{R}$ is linear or nonlinear depending on the linearity of $\mathscr{V}$. We focus below on the interface problem for the linear potential $\mathscr{V}=V(t,x)$, especially for $\mathscr{V}=V(x)$. 

For the transmission conditions presented in Section \ref{Sec_SWR}, we are going to show that if $\mathscr{V}=V(t,x)$, then
\begin{equation}
  \label{chp1_interfaceproblem}
  g^{k+1} = \mathcal{R} g^{k} = \mathcal{L} g^{k} + d,
\end{equation}
where $\mathcal{L}$ is a block matrix
\begin{equation} \label{chp1_L_N}
  \mathcal{L} = 
  \begin{pmatrix}
    & X^{2,1} & X^{2,2} & & & \\
    X^{1,4} \\
    & & & X^{3,1} & X^{3,2} \\
    & X^{2,3} & X^{2,4} \\
    & & & & & \cdots \\
    & & & X^{3,3} & X^{3,4} \\
    & & & & & & X^{N-1,1} & X^{N-1,2}\\
    & & & & &\cdots \\
    & & & & & & & & X^{N,1}\\
    & & & & & & X^{N-1,3} & X^{N-1,4}
  \end{pmatrix},
\end{equation}
with $X^{j,p} \in \mathbb{C}^{N_T \times N_T}$, $j=1,2,...,N$, $p=1,2,3,4$ and
$d$ is a vector
\begin{equation}
  \label{chp1_bN}
  d = \big( d_{1,r}^T, d_{2,l}^T, d_{2,r}^T, \cdots, d_{N,l}^T \big)^T
  \in \mathbb{C}^{(2N-2) \times N_T},
  \ 
  d_{j,l}, \ d_{j,r} \in \mathbb{C}^{N_T}.
\end{equation}
It is easy to see that the formula (\ref{chp1_interfaceproblem}) is equivalent to
 \begin{enumerate}
 \item for $j=1$,
   \begin{displaymath}
     \begin{pmatrix}
       l_{2,1}^{k+1}\\
       l_{2,2}^{k+1}\\
       \vdots\\
       l_{2,N_T}^{k+1}
     \end{pmatrix}
     = X^{1,4}
     \begin{pmatrix}
       r_{1,1}^{k}\\
       r_{1,2}^{k}\\
       \vdots\\
       r_{1,N_T}^{k}            
     \end{pmatrix}
     + d_{2,l},
   \end{displaymath}
 \item for $j=2,...,N-1$,
   \begin{equation}
     \label{chp1_Rlgb_lr}
     \begin{split}
       & \begin{pmatrix}
         r_{j-1,1}^{k+1}\\
         r_{j-1,2}^{k+1}\\
         \vdots\\
         r_{j-1,N_T}^{k+1}
       \end{pmatrix}
       = X^{j,1} 
       \begin{pmatrix}
         l_{j,1}^{k}\\
         l_{j,2}^{k}\\
         \vdots\\
         l_{j,N_T}^{k}      
       \end{pmatrix}
       + X^{j,2}
       \begin{pmatrix}
         r_{j,1}^{k}\\
         r_{j,2}^{k}\\
         \vdots\\
         r_{j,N_T}^{k}            
       \end{pmatrix}
       + d_{j-1,r},
       \\
       & \begin{pmatrix}
         l_{j+1,1}^{k+1}\\
         l_{j+1,2}^{k+1}\\
         \vdots\\
         l_{j+1,N_T}^{k+1}
       \end{pmatrix}
       = X^{j,3} 
       \begin{pmatrix}
         l_{j,1}^{k}\\
         l_{j,2}^{k}\\
         \vdots\\
         l_{j,N_T}^{k}      
       \end{pmatrix}
       + X^{j,4}
       \begin{pmatrix}
         r_{j,1}^{k}\\
         r_{j,2}^{k}\\
         \vdots\\
         r_{j,N_T}^{k}            
       \end{pmatrix}
       + d_{j+1,l},
     \end{split}
   \end{equation}
 \item for $j=N$,
   \begin{displaymath}
     \begin{pmatrix}
       r_{N-1,1}^{k+1}\\
       r_{N-1,2}^{k+1}\\
       \vdots\\
       r_{N-1,N_T}^{k+1}
     \end{pmatrix}
     = X^{N,1} 
     \begin{pmatrix}
       l_{N,1}^{k}\\
       l_{N,2}^{k}\\
       \vdots\\
       l_{N,N_T}^{k}      
     \end{pmatrix}
     + d_{N-1,r}.
   \end{displaymath}
 \end{enumerate}

\begin{proposition} 
  \label{chp1_PropN_S02}
%
For the transmission condition involving the operator $S_0^2$, in the case of linear potential $\mathscr{V}=V(t,x)$, if we assume that the matrices $\mathbb{A}_{j,n}- \mathbb{B}_{j,n}$  $n=1,2,...,N_T$ are not singular, then the $N$ equations \eqref{AvML} could be written in the global form of interface problem (\ref{chp1_interfaceproblem})
  $$
  g^{k+1} = \mathcal{L} g^{k} + d.$$
\end{proposition} 

\begin{pf} 
  First, according to \eqref{AvML}, we have
  \begin{align*}
    & (\mathbb{A}_{j,1} - \mathbb{B}_{j,1}) \mathbf{v}_{j,1}^{k}  = \frac{2i}{\Delta t} \mathbb{M}_j \mathbf{u}_{j,0} + e^{-\frac{i \pi}{4}} Q_j^T \sqrt{\frac{2}{\Delta t}}  \beta_1  Q_j \mathbf{v}_{j,0}^k - Q_j^T  (l_{j,1}^{k},r_{j,1}^{k})^T,\\
    & (\mathbb{A}_{j,n}  - \mathbb{B}_{j,n}) \mathbf{v}_{j,n}^{k}  = \frac{2i}{\Delta t} \mathbb{M}_j \mathbf{u}_{j,n-1}^k + e^{-\frac{i \pi}{4}}\sqrt{\frac{2}{\Delta t}} Q_j^T \sum_{q=0}^{n-1} \beta_{2-q} Q_j \mathbf{v}_{j,q}^k - Q_j^T (l_{j,n}^{k},r_{j,n}^{k})^T\\
    & \quad =  \frac{4i}{\Delta t} \mathbb{M}_j \mathbf{v}_{j,n-1}^k - \frac{2i}{\Delta t} \mathbb{M}_j \mathbf{u}_{j,n-2}^k 
    + e^{-\frac{i \pi}{4}}\sqrt{\frac{2}{\Delta t}} Q_j^T \sum_{q=0}^{n-1} \beta_{2-q} Q_j \mathbf{v}_{j,q}^k - Q_j^T (l_{j,n}^{k},r_{j,n}^{k})^T,\\ 
    & \quad = \sum_{q=1}^{n-1} \Big( (-1)^{n-1-q} \frac{4i}{\Delta t} \mathbb{M}_j + e^{-i\pi/4}\sqrt{\frac{2}{\Delta t}} \beta_{n-q} Q_j^T Q_j \Big) \mathbf{v}_{j,q}^k , \\  
	& \quad \quad +\Big( (-1)^{n-1} \frac{2i}{\Delta t}   \mathbb{M}_j + e^{-i\pi/4}\sqrt{\frac{2}{\Delta t}} \beta_{n} Q_j^T Q_j  \Big)\mathbf{u}_{j,0}- Q_j^T (l_{j,n}^{k},r_{j,n}^{k})^T, \ n \geqslant 2,
	\end{align*} 
  where we recall that $\mathbf{v}_{j,0}^k = \mathbf{u}_{j,0}$. Thus, we could see that  
  \begin{align}
    \label{vnvn_1}
    \mathbf{v}_{j,n}^k =&  -(\mathbb{A}_{j,n}  - \mathbb{B}_{j,n})^{-1} Q_j^T (l_{j,n}^{k},r_{j,n}^{k})^T \nonumber \\
    & + (\mathbb{A}_{j,n}  - \mathbb{B}_{j,n})^{-1} \sum_{q=1}^{n-1} \Big( (-1)^{n-1-q} \frac{4i}{\Delta t} \mathbb{M}_j + e^{-i\pi/4}\sqrt{\frac{2}{\Delta t}} \beta_{n-q} Q_j^T Q_j  \Big)  \mathbf{v}_{j,q}^k\\ 
    & + (\mathbb{A}_{j,n}  - \mathbb{B}_{j,n})^{-1} \Big( (-1)^{n-1} \frac{2i}{\Delta t}   \mathbb{M}_j + e^{-i\pi/4}\sqrt{\frac{2}{\Delta t}} \beta_{n} Q_j^T Q_j  \Big)\mathbf{u}_{j,0}.  \nonumber
  \end{align}  

  By induction on $n$, it is easy to see that $\mathbf{v}_{j,n}^k$ is a linear function of $l_{j,s}^{k}$ and $r_{j,s}^{k},s=1,2,...,n$. Then considering the formulas \eqref{TransmissionCond_disc}, in order to finish the proof, we need only verify that $\widetilde{S} (Q_{j,l} \mathbf{v}_{j,n}^{k})$ and $ \widetilde{S}(Q_{j,r} \mathbf{v}_{j,n}^{k})$ are linear functions of $\mathbf{v}_{j,s}^{k},s=1,2,...,n$. 
\end{pf}

\vbox{}
\begin{proposition} 
  \label{chp1_PropN_S}
  For any transmission condition presented in Section \ref{Sec_SWR}, assuming that the matrices $\mathbb{A}_{j,n}- \mathbb{B}_{j,n}$, $n=1,2,...,N_T$ are not singular, then the interface problem in the case of linear potential $\mathscr{V}=V(t,x)$ could be written in the global form (\ref{chp1_interfaceproblem}).
\end{proposition}
\begin{pf}
The proof is quite similar than that of the previous proposition. For each transmission condition, we only need to recalculate the expression of $\mathbf{v}_{j,n}^k$.
\end{pf}

We now turn to the structure of sub-blocks for $\mathscr{V}=V(x)$ and $j=2,3,...,N-1,$
\begin{align*}
  & X^{j,1} = \{ x^{j,1}_{n,s} \}_{1 \leqslant n,s \leqslant N_T}, \quad X^{j,2}=\{ x^{j,2}_{n,s} \}_{1 \leqslant n,s \leqslant N_T},\\
  & X^{j,3} = \{ x^{j,3}_{n,s} \}_{1 \leqslant n,s \leqslant N_T}, \quad X^{j,4}=\{ x^{j,4}_{n,s} \}_{1 \leqslant n,s \leqslant N_T}.
\end{align*}
and $X^{1,4} = \{ x^{1,4}_{n,s} \}_{1 \leqslant n,s \leqslant N_T}$ and $X^{N-1,1} = \{ x^{N-1,1}_{n,s} \}_{1 \leqslant n,s \leqslant N_T}$.
For 5 time steps, this structure is described below
\begin{displaymath}
  \begin{pmatrix}
    {\color{red}{\star}}\\
    {\color{blue}{\times}} & {\color{red}{\star}} \\
    {\color{green}{\circ}} & {\color{blue}{\times}} & {\color{red}{\star}} \\
    {\color{black}{\triangleleft}} & {\color{green}{\circ}} & {\color{blue}{\times}} & {\color{red}{\star}} \\
    {\color{violet}{\diamond}} & {\color{black}{\triangleleft}} & {\color{green}{\circ}} & {\color{blue}{\times}} & {\color{red}{\star}}
  \end{pmatrix}, \ N_T=5.
\end{displaymath}
thus, each sub-diagonal have an identical element.

\begin{proposition}
  \label{Prop_X_S02}
For the transmission condition involving the operator $S_0^2$, if $\mathscr{V} = V(x)$ and assuming that $\mathbb{A}_{j,n}- \mathbb{B}_{j,n}$, $n=1,2,...,N_T$ are not singular, then the matrices $X^{1,4}$ $X^{j,1},X^{j,2},X^{j,3},X^{j,4},j=2,3,...,N-1$ and $X^{N,1}$ are lower triangular matrices and they satisfy
  \begin{displaymath}
    \begin{split}
      & x^{1,4}_{n,s} = x^{1,4}_{n-1,s-1},\\
      & x^{j,1}_{n,s} = x^{j,1}_{n-1,s-1},\
      x^{j,2}_{n,s} = x^{j,2}_{n-1,s-1},\\
      & x^{j,3}_{n,s} = x^{j,3}_{n-1,s-1},\
      x^{j,4}_{n,s} = x^{j,4}_{n-1,s-1}, j=2,3,...,N-1,\\
      & x^{N,1}_{n,s} = x^{N,1}_{n-1,s-1},
    \end{split}
  \end{displaymath} 
  for $2 \leqslant s \leqslant n \leqslant N_T $.
\end{proposition}

\begin{pf}
Without loss of generality, we consider here $ j = 2,3, ..., N-1 $. First, we design
  \begin{displaymath}
    \mathbb{Y}_{n,q}^j = 
    \left \{ 
      \begin{array}{ll}
        -(\mathbb{A}_{j,n}  - \mathbb{B}_{j,n})^{-1}, \ q=n,\\ 
        (\mathbb{A}_{j,n}  - \mathbb{B}_{j,n})^{-1} \Big( (-1)^{n-1-q} \frac{4i}{\Delta t} \mathbb{M}_j + e^{\frac{-i\pi}{4}}\sqrt{\frac{2}{\Delta t}} \beta_{n-q} Q_j^T Q_j  \Big), \ q=1,2,...,n-1.
      \end{array}
    \right.
  \end{displaymath}
If the linear potential $\mathscr{V} = V(x)$ is independent of time, then it is easy to see
  \begin{displaymath}
    \mathbb{A}_{j,1} = \mathbb{A}_{j,2} = \cdots = \mathbb{A}_{j,N_T},\
    \mathbb{B}_{j,1} = \mathbb{B}_{j,2} = \cdots = \mathbb{B}_{j,N_T}.
  \end{displaymath}
  Thus for $2 \leqslant s \leqslant n \leqslant N_T $,
  \begin{equation}
    \label{chp1_Ynn_1}
    \mathbb{Y}_{n,s}^j = \mathbb{Y}_{n-1,s-1}^j.
  \end{equation}

  Then, according to \eqref{vnvn_1}, we have
  \begin{equation}
    \label{vYU}
    \mathbf{v}_{j,n}^k = \mathbb{Y}_{n,n}^j Q_j^T (l_{j,n}^{k},r_{j,n}^{k})^T + \sum_{q=1}^{n-1} \mathbb{Y}_{n,q}^j \mathbf{v}_{j,q}^k +  \mathbb{U}_{j,n} \mathbf{u}_{j,0}.
  \end{equation}
  where $\mathbb{U}_{j,n}=(\mathbb{A}_{j,n}  - \mathbb{B}_{j,n})^{-1} \Big( \frac{2i}{\Delta t}  (-1)^{n-1} \mathbb{M}_j + e^{-i\pi/4}\sqrt{\frac{2}{\Delta t}} \beta_{n} Q_j^T Q_j  \Big)$. By induction, we can obtain an expression of $\mathbf{v}_{j,n}^k$:
  \begin{equation}
    \label{vLU}
    \mathbf{v}_{j,n}^k = \sum_{q=1}^n \mathbb{L}_{n,q}^j Q_j^T (l_{j,q}^k,r_{j,q}^k)^T + U_{j,n} \mathbf{u}_{j,0},
  \end{equation}
  where $\mathbb{L}_{n,q}^j$, $q=1,2,...,n$ and $U_{j,n}$ are matrix. For example, $\mathbb{L}_{n,n}^j = \mathbb{Y}_{n,n}^j$. We are going to show that for $2 \leqslant s \leqslant n \leqslant N_T$,
  \begin{equation}
    \label{LLn_1}
    \mathbb{L}_{n,s}^j = \mathbb{L}_{n-1,s-1}^j.
  \end{equation}
  Replacing $\mathbf{v}_{j,q}^k$ in \eqref{vYU} by \eqref{vLU}, we have
  \begin{align*}
    \mathbf{v}_{j,n}^k =& \mathbb{Y}_{n,n}^j Q_j^T (l_{j,n}^{k},r_{j,n}^{k})^T + \sum_{q=1}^{n-1} \mathbb{Y}_{n,q}^j \Big( \sum_{p=1}^q \mathbb{L}_{q,p}^j Q_j^T (l_{j,p}^k,n_{j,p}^k)^T + U_{j,q} \mathbf{u}_{j,0} \Big)  +  \mathbb{U}_{j,n} \mathbf{u}_{j,0}\\
    = & \mathbb{Y}_{n,n}^j Q_j^T (l_{j,n}^{k},r_{j,n}^{k})^T + \sum_{p=1}^{n-1} \Big( \sum_{q=p}^{n-1} \mathbb{Y}_{n,q}^j \mathbb{L}_{q,p}^j \Big) Q_j^T (l_{j,p}^k,r_{j,p}^k)^T + \Big( \sum_{q=1}^{n-1} \mathbb{Y}_{n,q}^j U_{j,q} +  \mathbb{U}_{j,n} \Big) \mathbf{u}_{j,0}.
  \end{align*}
  Comparing the above formula with \eqref{vLU}, we have
  \begin{equation}
    \mathbb{L}_{n,s} = 
    \left \{
      \begin{array}{ll}
        \begin{aligned}
          \mathbb{Y}_{n,n}^j,
        \end{aligned} \\
        \begin{aligned}
          \sum_{q=s}^{n-1} \mathbb{Y}_{n,q}^j \mathbb{L}_{q,s}^j, 1 \leqslant s < n,
        \end{aligned}
      \end{array}
    \right.
    \Rightarrow
    \mathbb{L}_{n-1,s-1} = 
    \left \{
      \begin{array}{ll}
        \begin{aligned}
          \mathbb{Y}_{n-1,n-1}^j,
        \end{aligned} \\
        \begin{aligned}
          \sum_{q=s-1}^{n-2} \mathbb{Y}_{n-1,q}^j \mathbb{L}_{q,s-1}^j, 2 \leqslant s < n.
        \end{aligned}
      \end{array}
    \right.
  \end{equation}
  By using \eqref{chp1_Ynn_1} and by induction on $n$, we get
  \begin{gather*}
    \mathbb{L}_{n,n} = \mathbb{Y}_{n,n}  = \mathbb{Y}_{n-1,n-1} =  \mathbb{L}_{n-1,n-1},\\
    \mathbb{L}_{n,s} = \sum_{q=s}^{n-1} \mathbb{Y}_{n,q}^j \mathbb{L}_{q,s}^j
    = \sum_{q=s}^{n-1} \mathbb{Y}_{n-1,q-1}^j \mathbb{L}_{q-1,s-1}^j
    = \sum_{q=s-1}^{n-2} \mathbb{Y}_{n-1,q}^j \mathbb{L}_{q,s-1}^j = \mathbb{L}_{n-1,s-1}, 2 \leqslant s < n. 
  \end{gather*}
  The formula \eqref{LLn_1} is thus demonstrated.

Then we replace $\mathbf{v}_{j,k}^n $ in the first two formulas of \eqref{TransmissionCond_disc} by \eqref{vLU}. We get
  \begin{align}
    \label{chp1_lrL}
    l_{j+1,n}^{k+1} =&  -r_{j,n}^{k} + 2 e^{-i\pi/4}\sqrt{\frac{2}{\Delta t}} \sum_{p=1}^{n}\beta_{n-p} \sum_{q=1}^{p} \mathbb{L}_{p,q}^j  Q_j^T (l_{j,q}^k,r_{j,q}^k)^T + R_{l,j,n}^k \nonumber \\
    = & -r_{j,n}^{k} + 2 e^{-i\pi/4}\sqrt{\frac{2}{\Delta t}} \sum_{q=1}^{n} Q_{j,r} \Big( \sum_{p=q}^{n} \beta_{n-p} \mathbb{L}_{p,q}^j  \Big)Q_j^T  (l_{j,q}^k,r_{j,q}^k)^T + R_{l,j,n}^k,  \\
    r_{j-1,n}^{k+1} =& -l_{j,n}^{k} + 2 e^{-i\pi/4}\sqrt{\frac{2}{\Delta t}} \sum_{q=1}^{n} Q_{j,l} \Big( \sum_{p=q}^{n} \beta_{n-p} \mathbb{L}_{p,q}^j  \Big) Q_j^T (l_{j,q}^k,r_{j,q}^k)^T + R_{r,j,n}^k, \nonumber
  \end{align}
  where we denotes the terms that are independent of $l_{j,s}^k$ and $r_{j,s}^k$, $s=1,2,...,N_T$ by remainder terms $R_{l,r}$ to make the proof more readable.

  Moreover, according to (\ref{chp1_Rlgb_lr}), we have 
  \begin{gather*}
    l_{2,n}^{k+1} = \sum_{s=1}^{N_T} x^{1,4}_{n,s} r_{1,s}^k + d_{2,l,n}, \ r_{N-1,n}^{k+1} = \sum_{s=1}^{N_T} x^{N,1}_{n,s} l_{N,s}^k + d_{N-1,r,n}, \\
    r_{j-1,n}^{k+1} = \sum_{s=1}^{N_T} x^{j,1}_{n,s} l_{j,s}^k + \sum_{s=1}^{N_T} x^{j,2}_{n,s} r_{j,s}^k + d_{j-1,r,n}, \\ 
    l_{j+1,n}^{k+1} = \sum_{s=1}^{N_T} x^{j,3}_{n,s} l_{j,s}^k + \sum_{s=1}^{N_T} x^{j,4}_{n,s} r_{j,s}^k + d_{j+1,l,n}.
  \end{gather*}
where $d_{j-1,l,n}$ and $d_{j+1,r,n}$ denote the $n$-th element of $d_{j-1,l}$ and $d_{j+1,r}$ respectively.

  Comparing the above formula with \eqref{chp1_lrL}, we have for $1 \leqslant n<s \leqslant N_T$,
  \begin{displaymath}
    x^{j,1}_{n,s} = x^{j,2}_{n,s} = x^{j,3}_{n,s} = x^{j,4}_{n,s} = 0,
  \end{displaymath}
  and for $1 \leqslant s \leqslant n \leqslant N_T$,
  \begin{align*}
    & x^{j,1}_{n,s} = -1 + 2 c_2 Q_{j,l} \Big( \sum_{p=s}^{n} \beta_{n-p} \mathbb{L}_{p,s}^j  \Big) Q_{j,l}^T,\ x^{j,2}_{n,s} = 2 c_2 Q_{j,l} \Big( \sum_{p=s}^{n} \beta_{n-p} \mathbb{L}_{p,s}^j  \Big) Q_{j,r}^T,\\
    & x^{j,3}_{n,s} = 2 c_2 Q_{j,r} \Big( \sum_{p=s}^{n} \beta_{n-p} \mathbb{L}_{p,s}^j  \Big) Q_{j,l}^T,\
    x^{j,4}_{n,s} = -1 + 2 c_2 Q_{j,r} \Big( \sum_{p=s}^{n} \beta_{n-p} \mathbb{L}_{p,s}^j  \Big) Q_{j,r}^T,
  \end{align*}
  where $c_2=e^{-i\pi/4}\sqrt{\frac{2}{\Delta t}}$ and we use $ Q_{j}^T (l_{j,q}^k,r_{j,q}^k)^T = Q_{j,l}^T l_{j,q}^k + Q_{j,r}^T r_{j,q}^k$. 

  Finally, using \eqref{LLn_1}, we have for $2 \leqslant s \leqslant n \leqslant N_T$,
  \begin{align*}
    x^{j,1}_{n,s} & = -1 + 2 e^{-i\pi/4}\sqrt{\frac{2}{\Delta t}}  Q_{j,l} \Big( \sum_{p=s}^{n} \beta_{n-p} \mathbb{L}_{p,s}^j  \Big) Q_{j,l}^T \\
    & = -1 + 2 e^{-i\pi/4}\sqrt{\frac{2}{\Delta t}} \ Q_{j,l} \Big( \sum_{p=s}^{n} \beta_{n-p} \mathbb{L}_{p-1,s-1}^j  \Big) Q_{j,l}^T \\
    & = -1 + 2 e^{-i\pi/4}\sqrt{\frac{2}{\Delta t}} Q_{j,l} \Big( \sum_{p=s-1}^{n-1} \beta_{n-1-p} \mathbb{L}_{p,s}^j  \Big) Q_{j,l}^T = x^{j,1}_{n-1,s-1}.
  \end{align*}
  In the same way, we can prove that $x^{j,2}_{n,s} = x^{j,2}_{n-1,s-1}$, $x^{j,3}_{n,s} = x^{j,3}_{n-1,s-1}$ and $x^{j,4}_{n,s} = x^{j,4}_{n-1,s-1}$.
\end{pf}

\begin{proposition}
  \label{Prop_X_S}
	With any transmission condition presented in Section \ref{Sec_SWR}, if  $\mathscr{V} = V(x)$ and assuming that $\mathbb{A}_{j,n}- \mathbb{B}_{j,n}$, $n=1,2,...,N_T$ are not singular, then the matrices $X^{1,4}$, $X^{j,1}$, $X^{j,2}$, $X^{j,3}$, $X^{j,4}$, $j=2,3,...,N-1$ and $X^{N,1}$ are lower triangular matrices and they satisfy
  \begin{displaymath}
    \begin{split}
      & x^{1,4}_{n,s} = x^{1,4}_{n-1,s-1},\\
      & x^{j,1}_{n,s} = x^{j,1}_{n-1,s-1},\
      x^{j,2}_{n,s} = x^{j,2}_{n-1,s-1},\\
      & x^{j,3}_{n,s} = x^{j,3}_{n-1,s-1},\
      x^{j,4}_{n,s} = x^{j,4}_{n-1,s-1}, j=2,3,...,N-1,\\
      & x^{N,1}_{n,s} = x^{N,1}_{n-1,s-1},
    \end{split}
  \end{displaymath} 
  for $2 \leqslant s \leqslant n \leqslant N_T $.
\end{proposition}
%
\begin{pf}
{The proof is similar to that of Proposition \ref{Prop_X_S02}. We only need to recompute $\mathbf{v}_{j,n}^k$ and $\mathbb{Y}^j_{n,q}$ for each transmission condition.}

\end{pf}

\section{New algorithm for time independent linear potential}

The standard implementation of the SWR method for the time-independent equations leads to the following classical algorithm

\begin{algorithm}[H]
  \caption{Classical algorithm}
  \begin{algorithmic}[1] 
    \STATE Initialize the iteration by $g^0$, 
    \STATE Solve Schr{\"o}dinger on each subdomain with $g^k$.
    \STATE Exchange values at interfaces and compute $g^{k+1}$.
    \STATE Do again steps 2 and 3 until error $||g^{k+1} - g^{k}||<\varepsilon$, $\varepsilon \ll 1$.
  \end{algorithmic}
\end{algorithm}

\vspace*{2mm}

As we can see, the classical algorithm requires to solve $K$ times the Schr{\"o}dinger equation on each subdomain, where $K$ corresponds to the number of iterations required to reach convergence. We are going to present a new algorithm for $\mathscr{V}=V(x)$ which is more efficient. As we will see,  it will require to solve the Schr{\"o}dinger equation on each subdomain only four times in total. This new algorithm is equivalent to the classical algorithm, but it reduces significantly the calculations.

Before giving this new algorithm, we could see that the classical algorithm is based on \eqref{interfacepbR}: $g^{k+1}=\mathcal{R} g^{k}$,
where the operator $\mathcal{R}$ includes the steps 2 and 3. We have shown in Proposition \ref{chp1_PropN_S} that
\begin{equation}
  \label{chp1_RLd}
  g^{k+1}=\mathcal{R} g^{k} = \mathcal{L} g^{k} + d.
\end{equation}
It is easy to see that \eqref{chp1_RLd} is nothing but the fix point method to solve the equation
\begin{equation}
  \label{chp1_ILgd}
  (I - \mathcal{L}) g = d.
\end{equation}
A big advantage to interpret \eqref{chp1_RLd} as a fixed point method to solve\eqref{chp1_ILgd} is that we can use any other iterative methods to solve this linear system. So we can use Krylov methods (ex. Gmres, Bicgstab) \cite{Saad2003}, which could accelerate the convergence prospectively. To use the Krylov methods or fixed point method, it is enough to define the application of $I - \mathcal{L}$ to vector $g$ by
\begin{displaymath}
  (I - \mathcal{L}) g = I - \mathcal{R}g + d.
\end{displaymath}
The classical algorithm could then be rewritten with

\begin{algorithm}[H]
  \caption{Classical algorithm, version 2}
  \begin{algorithmic}[1] 
    \STATE Build $d = \mathcal{R} \cdot \mathbf{0}$ in (\ref{chp1_ILgd}) explicitly,  
    \STATE Define the application of $I - \mathcal{L}$ to vector in (\ref{chp1_ILgd}),
    \STATE Solve the linear system (\ref{chp1_ILgd}) by an iterative method (fixed point or Krylov).
    \STATE Solve the Schr{\"o}dinger equation on each subdomain for each time step using the boundary conditions obtained at step 3.
  \end{algorithmic}
\end{algorithm}

\vspace*{2mm}
If the fixed point method is used in Step 3, we recover the first version of the classical algorithm. The second version of the classical algorithm allows the use of Krylov methods to accelerate convergence. However, applying $(I-\mathcal{L})$ to vector $g$ is still a very expensive operation. With the help of Propositions \ref{Prop_X_S02} and \ref{Prop_X_S}, we propose a new algorithm 

\begin{algorithm}[H]
  \caption{New algorithm}
  \label{chp1_algo_new}
  \begin{algorithmic}[1] 
    \STATE Build $\mathcal{L}$ and $d$ in (\ref{chp1_ILgd}) explicitly, 
    \STATE Solve (\ref{chp1_ILgd}) by an iterative method,
    \STATE Solve Schr{\"o}dinger equation on each subdomain using the boundary conditions obtained at step 2.
  \end{algorithmic}
\end{algorithm}

\vspace*{2mm}

We show beloa the construction of the matrix $\mathcal{L}$ and the vector $d$. As it will be seenn, their computation is not costly. Regarding the implementation, we then show how $\mathcal{L} $ and $d$ are stored for use of parallelism. Here, we use the PETSc library \cite{petsc-user-ref}. Using the matrix form in PETSc, the memory required for each MPI process \cite{mpiforum30} is independent of the number of subdomains.

\subsection{Construction of the matrix $\mathcal{L}$ and the vector $d$}

We use the formulas (\ref{AvML}) and (\ref{TransmissionCond_disc}) for the constructions. Numerically, we consider $l_{j,n}^k$ and $r_{j,n}^k$ as inputs, and $l_{j-1,n}^{k+1}$ and $r_{j+1,n}^{k+1}$ as outputs:
\begin{displaymath}
  \text{inputs: } l_{j,n}^k, r_{j,n}^k \longrightarrow \text{(\ref{TransmissionCond_disc})} \longrightarrow \text{outputs: } l_{j-1,n}^{k+1}, r_{j+1,n}^{k+1}. 
\end{displaymath}

It is easy to see that
\begin{displaymath}
  d = \big( d_{1,r}^T, d_{2,l}^T, d_{2,r}^T, \cdots, d_{N,l}^T \big)^T
  = \mathcal{R} \cdot \bf{0},
\end{displaymath}
where $\bf{0}$ is the zero vector. The elements of $d$ are obtained by
\begin{displaymath}
  d_{j-1,r} = 
  \begin{pmatrix}
    r_{j-1,1}^{k+1}\\
    r_{j-1,2}^{k+1}\\
    \vdots\\
    r_{j-1,N_T}^{k+1}
  \end{pmatrix},\quad
  d_{j+1,l} = 
  \begin{pmatrix}
    l_{j+1,1}^{k+1}\\
    l_{j+1,2}^{k+1}\\
    \vdots\\
    l_{j+1,N_T}^{k+1}
  \end{pmatrix},
\end{displaymath}
where the scalars $r_{j-1,s}^{k+1}, l_{j+1,s}^{k+1},s=1,2,...,N_T$ are given by the formula (\ref{TransmissionCond_disc}) with
\begin{displaymath}
  l^k_{j,s}=r^k_{j,s}=0, s=1,2,...,N_T.
\end{displaymath}
The equation is solved numerically on each subdomain only one time. Note that this construction works for $\mathscr{V} = V(t,x)$.

According to Propositions \ref{Prop_X_S} and \ref{Prop_X_S02}, if $\mathscr{V}=V(x)$, in order to build the matrix $\mathcal{L}$, it is enough to compute the first columns of blocks $X^{1,4}$, $X^{j,1}$, $X^{j,2}$, $X^{j,3}$, $X^{j,4}$, $j=2,3,...,N-1$ and $X^{N,1}$.

The first column of $X^{j,1}$ is
\begin{displaymath}
  X^{j,1}
  \begin{pmatrix}
    1 \\
    0 \\
    \vdots \\
    0
  \end{pmatrix}
  = 
  \Big(
  X^{j,1}
  \begin{pmatrix}
    1 \\
    0 \\
    \vdots \\
    0
  \end{pmatrix}
  + 
  X^{j,2}
  \begin{pmatrix}
    0 \\
    0 \\
    \vdots \\
    0
  \end{pmatrix}
  +
  d_{j-1,r} \Big)
  - d_{j-1,r}
  = 
  \begin{pmatrix}
    r_{j-1,1}^{k+1}\\
    r_{j-1,2}^{k+1}\\
    \vdots\\
    r_{j-1,N_T}^{k+1}
  \end{pmatrix}
  - d_{j-1,r}.
\end{displaymath}
The first column of $X^{j,3}$ is
\begin{displaymath}
  X^{j,3}
  \begin{pmatrix}
    1 \\
    0 \\
    \vdots \\
    0
  \end{pmatrix}
  = 
  \Big(
  X^{j,3}
  \begin{pmatrix}
    1 \\
    0 \\
    \vdots \\
    0
  \end{pmatrix}
  + 
  X^{j,4}
  \begin{pmatrix}
    0 \\
    0 \\
    \vdots \\
    0
  \end{pmatrix}
  +
  d_{j+1,l} \Big)
  - d_{j+1,l}
  = 
  \begin{pmatrix}
    l_{j+1,1}^{k+1}\\
    l_{j+1,2}^{k+1}\\
    \vdots\\
    l_{j+1,N_T}^{k+1}
  \end{pmatrix}
  - d_{j+1,l}.
\end{displaymath}
The scalars $r_{j-1,s}^{k+1}, l_{j+1,s}^{k+1},s=1,2,...,N_T$ are computed by the formula (\ref{TransmissionCond_disc}) with
\begin{displaymath}
  l^k_{j,s}=r^k_{j,s}=0, s=1,2,...,N_T\ \text{except for} \ l_{j,1}^k=1.
\end{displaymath}
The equation is solved numerically only one time on the subdomain $(a_j,b_j)$.

In the same way, the first columns of $X^{j,2}$ and $X^{j,4}$ are
\begin{displaymath}
  X^{j,2}
  \begin{pmatrix}
    1 \\
    0 \\
    \vdots \\
    0
  \end{pmatrix}
  = 
  \Big(
  X^{j,2}
  \begin{pmatrix}
    1 \\
    0 \\
    \vdots \\
    0
  \end{pmatrix}
  + 
  X^{j,4}
  \begin{pmatrix}
    0 \\
    0 \\
    \vdots \\
    0
  \end{pmatrix}
  +
  d_{j-1,r} \Big)
  - d_{j-1,r}
  = 
  \begin{pmatrix}
    r_{j-1,1}^{k+1}\\
    r_{j-1,2}^{k+1}\\
    \vdots\\
    r_{j-1,N_T}^{k+1}
  \end{pmatrix}
  - d_{j-1,r},
\end{displaymath}
and
\begin{displaymath}
  X^{j,4}
  \begin{pmatrix}
    1 \\
    0 \\
    \vdots \\
    0
  \end{pmatrix}
  = 
  \Big(
  X^{j,2}
  \begin{pmatrix}
    0 \\
    0 \\
    \vdots \\
    0
  \end{pmatrix}
  + 
  X^{j,4}
  \begin{pmatrix}
    1 \\
    0 \\
    \vdots \\
    0
  \end{pmatrix}
  +
  d_{j+1,l} \Big)
  - d_{j+1,l}
  = 
  \begin{pmatrix}
    l_{j+1,1}^{k+1}\\
    l_{j+1,2}^{k+1}\\
    \vdots\\
    l_{j+1,N_T}^{k+1}
  \end{pmatrix}
  - d_{j+1,l},
\end{displaymath}
where the scalars $r_{j-1,s}^{k+1}, l_{j+1,s}^{k+1},s=1,2,...,N_T$ are obtained by the formula (\ref{TransmissionCond_disc}), but with
\begin{displaymath}
  l^k_{j,s}=r^k_{j,s}=0, s=1,2,...,N_T\ \text{except for} \ r_{j,1}^k=1.
\end{displaymath}
The equation is solved numerically on each subdomain $(a_j,b_j)$ only one time.

In conclusion, it is sufficient to solve the equation {\eqref{Algo}} on each subdomain three times to construct explicitly the interface problem. The construction is inexpensive. In total, the equation {\eqref{Algo}} is solved on each subdomain four times in the new algorithm. Numerically, we will compare the classical and the new algorithms in Section \ref{Sec_compare_algorithmes}.

\subsection{Storage of the matrix $\mathcal{L}$ and the vector $d$ for massive parallel computing}
\label{chp1_chp_stokageLb}

Thanks to the peculiar form of the matrix $\mathcal{L}$, we can build it on parallel computers through an MPI implementation. The transpose of $\mathcal{L}$ is stored in a distributed manner using the library PETSc. As we can see below, the first block column of $\mathcal{L}$ is in MPI process 0. The second and third blocks columns are in MPI process 1, and so on for other processes. The consumed  memory for each process is at most the sum of 4 blocks. The size of each block is $N_T \times N_T$. Each block contain $(N_T+1) \times N_T/2$ non zero elements according to Propositions \ref{chp1_PropN_S} and \ref{chp1_PropN_S02}.

\begin{equation}
  \label{chp1_Lmpi}
  \mathcal{L} = 
  \begin{pmatrix}
    \multicolumn{1}{l}{\overbrace{\hspace{2.0em}}^{\mathrm{MPI}\ 0}} & 
    \multicolumn{2}{l}{\overbrace{\hspace{5.0em}}^{\mathrm{MPI}\ 1}} &
    \multicolumn{2}{l}{\overbrace{\hspace{5.0em}}^{\mathrm{MPI}\ 2}} &
    & 
    \multicolumn{2}{l}{\overbrace{\hspace{8.0em}}^{\mathrm{MPI}\ N-2}} &
    \multicolumn{2}{l}{\overbrace{\hspace{2.0em}}^{\mathrm{MPI}\ N-1}}
    \\
    & X^{2,1} & X^{2,2} & & & \\
    X^{1,4} \\
    & & & X^{3,1} & X^{3,2} \\
    & X^{2,3} & X^{2,4} \\
    & & & & & \cdots \\
    & & & X^{3,3} & X^{3,4} \\
    & & & & & & X^{N-1,1} & X^{N-1,2}\\
    & & & & &\cdots \\
    & & & & & & & & X^{N,1}\\
    & & & & & & X^{N-1,3} & X^{N-1,4}
  \end{pmatrix}.
\end{equation}

The vector $d$ can also be stored in PETSc form. The first block is in MPI process 0, the second and the third are in MPI process 1, and so on. The last block is in MPI process $N-1$. Each MPI process contain at most $2\times N_T$ elements.
\begin{displaymath}
  d = \big( \underbrace{d_{1,r}^T}_{\text{MPI 0}}, \underbrace{d_{2,l}^T, d_{2,r}^T}_{\text{MPI 1}}, \cdots, \underbrace{d_{j,l}^T, d_{j,r}^T}_{\text{MPI $j-1$}}, \cdots, \underbrace{d_{N,l}^T}_{\text{MPI $N-1$}} \big)^T.
\end{displaymath}

\section{Preconditioned algorithm for general potentials}

In Section \ref{Sec_InterfacePb}, we have established the interface problem \eqref{interfacepbR} for Schr{\"o}dinger equation with time dependent or nonlinear potential. However, it is not possible to construct the interface matrix $\mathcal{L}$ without much computation since the Propositions \ref{Prop_X_S02} and \ref{Prop_X_S} only hold for time independent linear potential. Thus, the new algorithm is not suitable here. Instead, to reduce the number of iterations required for convergence, we propose to add a preconditioner $P^{-1}$ ($P$ is a non singular matrix) in \eqref{interfacepbR} which leads to the preconditioned algorithm:
\begin{enumerate}
\item for $\mathscr{V}=V(t,x)$,
  \begin{align}
    & g^{k+1} = I - P^{-1} (I - \mathcal{R}) g^{k}, \label{chp2_algopd_L} \\
    & P^{-1} (I - \mathcal{L}) = P^{-1} d,   \label{chp2_algopd_Lpf} 
  \end{align}
\item for $\mathscr{V}=f(u)$,
  \begin{equation}
    \label{chp2_algopd_NL}
    g^{k+1} = I - P^{-1} (I - \mathcal{R}_{nl}) g^{k}.
  \end{equation}
  \end{enumerate}	

We now turn to explain which preconditioner is used.
The interface problem for the free Schr\"odinger equation (without potential) is 
\begin{displaymath}
  g^{k+1} = \mathcal{L}_0 g^{k} + d,
\end{displaymath}
where the symbol $\mathcal{L}_0$ is used to highlight here the potential is zero. The transmission condition is the same as that for \eqref{Schequ}. We propose for time dependent or nonlinear potential the preconditioner as
\begin{displaymath}
P = I - \mathcal{L}_0.
\end{displaymath}
We have two reasons to believe that this is a good choice.
%
\begin{enumerate}
\item
  The matrix $\mathcal{L}_0$ can be constructed easily since a zero potential is independent of time. Therefore, the construction of $\mathcal{L}_0$ only needs to solve the free Schr{\"o}dinger equation two times on each subdomains. This construction is therefore scalable.
\item
 Intuitively, the Schr{\"o}dinger operator without potential is a roughly approximating of the Schr{\"o}dinger operator with potential:
  \begin{displaymath}
    i\partial_t + \partial_{xx} \approx i\partial_t  + \partial_{xx}u  + \mathscr{V},
  \end{displaymath}
  thus
  \begin{displaymath}
      P = I - \mathcal{L}_0  \approx I - \mathcal{L},\quad
      P = I - \mathcal{L}_0  \approx I - (\mathcal{R}_{nl}-\mathcal{R}_{nl} \cdot \mathbf{0}).
  \end{displaymath}
\end{enumerate}

Next, we present the application of preconditioner. The transpose of $P$ is stored in PETSc form. For any vector $y$, the vector $x:=P^{-1} y$ is computed by solving the linear system
\begin{equation}
  \label{Pxg}
  P x = (I - \mathcal{L}_0) x= y \ \Leftrightarrow \ x^T P^T = y^T.
\end{equation}   
We do not explicitly construct the matrix $P^{-1}$ as the inverse of a distributed matrix numerically is too expensive. The linear system (\ref{Pxg}) is solved by the Krylov methods (Gmres or Bicgstab) initialized by zero vector using the library PETSc. We will see in Section \ref{Sec_compare_conditions} that the computation time for applying this preconditioner is quite small compared with the computation time for solving the Schr{\"o}dinger equation on subdomains. 

\section{Numerical results}

\begin{figure}[!htbp]
  \centering
  \includegraphics[width=.46\textwidth]{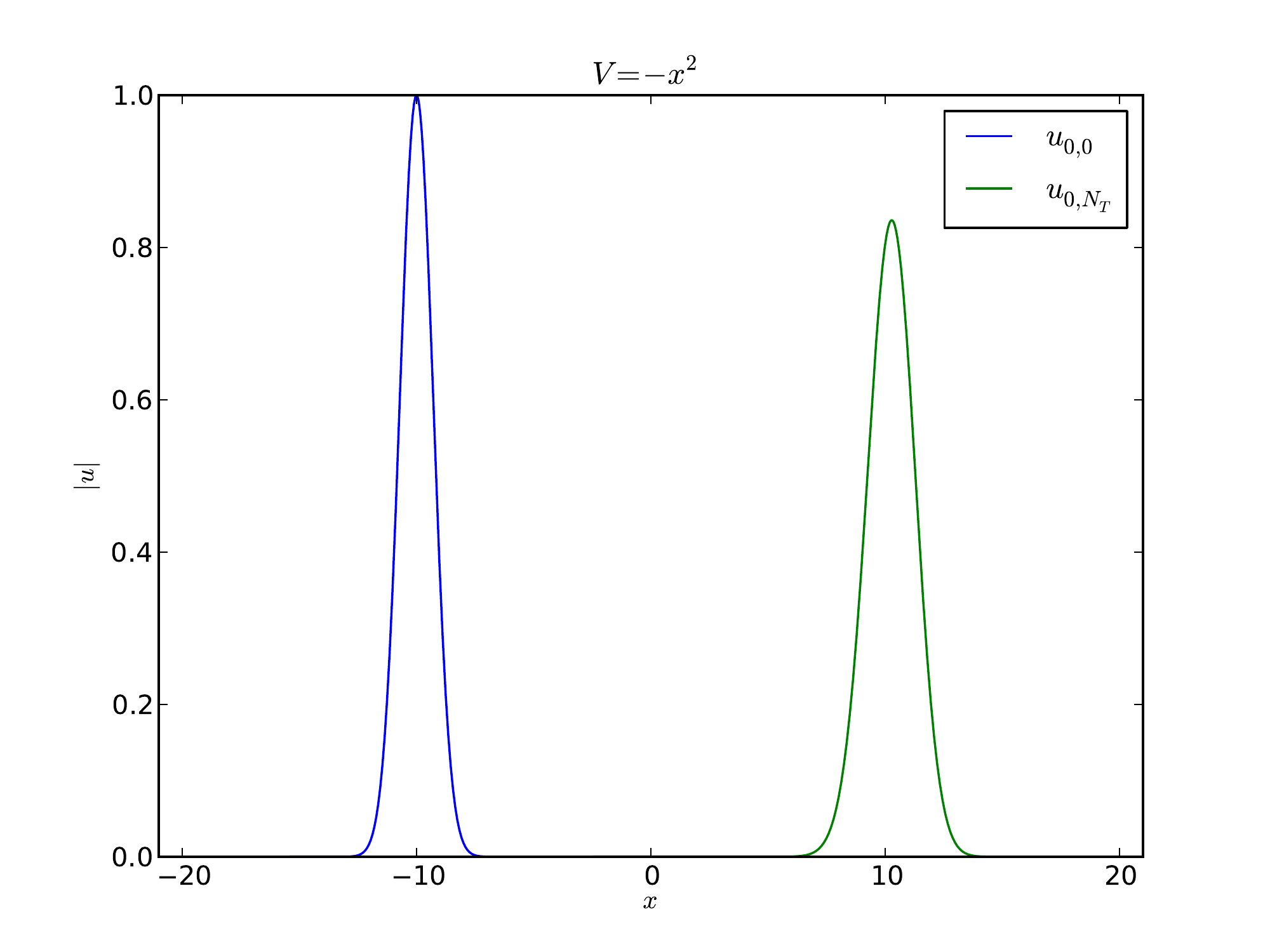}
  \includegraphics[width=.46\textwidth]{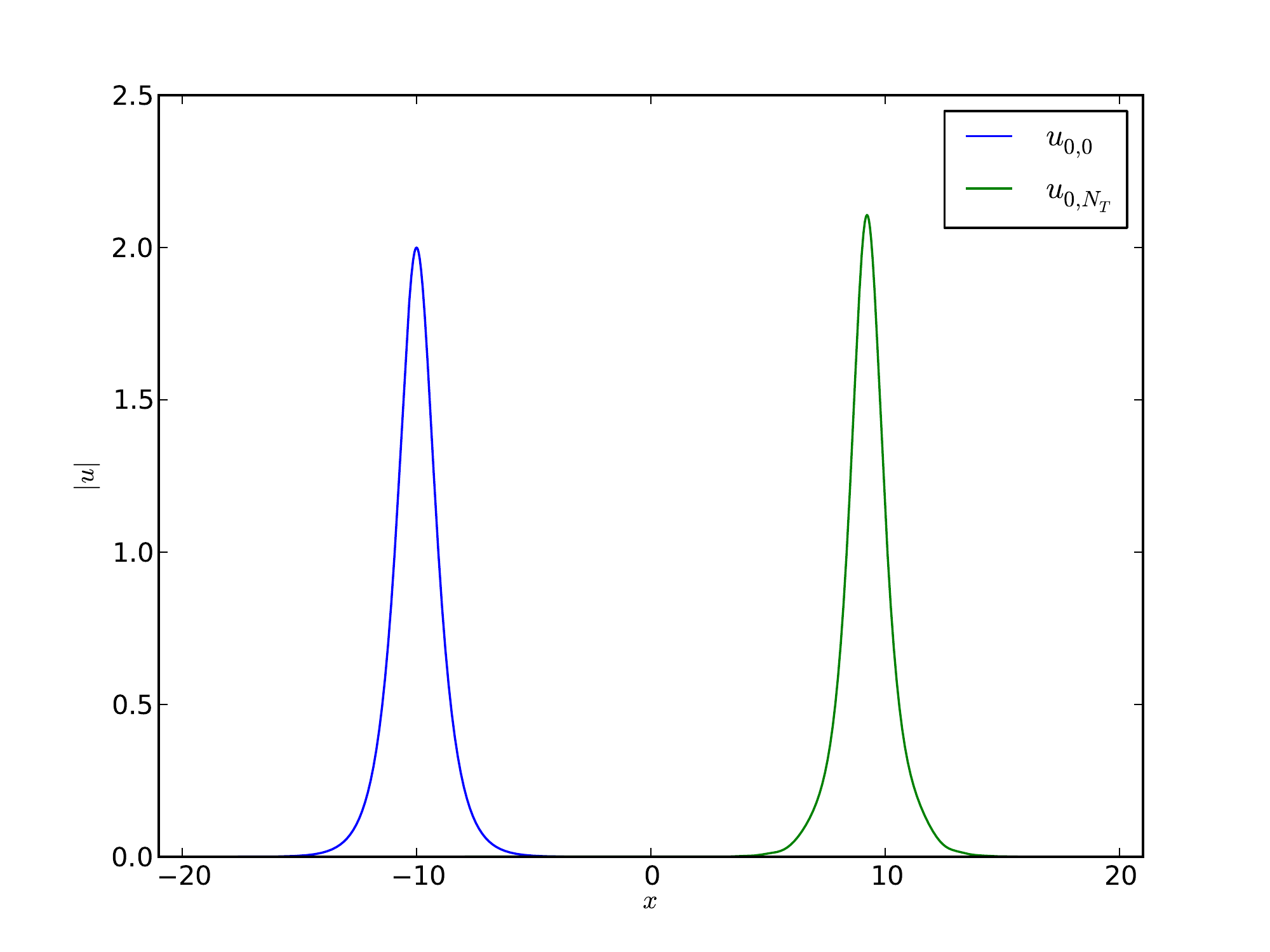}
  \caption{$|u_{0,0}|$ et $|u_{0,N_T}|$ on $(a_0, b_0)$, $\mathscr{V} = -x^2$ (left) and $\mathscr{V} = |u|^2$ (right), $\Delta t=0.001$, $\Delta x = 10^{-5}$.}
  \label{Initialsol}
\end{figure}

The physical domain $(a_0,b_0)=(-21,21)$ is decomposed into $N$ equal subdomains without overlap. We fix in this section the final time to $T=0.5$, the time step to $\Delta t=0.001$ and the mesh size to $\Delta x = 10^{-5}$ without special statement. The potentials that we consider in this part and the corresponding initial data are
\begin{enumerate}
\item time independent linear potential: $\mathscr{V}=-x^2$, $u_0(x) = e^{-(x+10)^2+20i(x+10)}$,
\item time dependent linear potential: $\mathscr{V}=5tx$, $u_0(x) = e^{-(x+10)^2+20i(x+10)}$,
\item nonlinear potential: $\mathscr{V}=|u|^2$, $u_0(x) = 2 \mathrm{sech}\big(\sqrt{2}(x+10)\big)e^{20i (x+10)}$,
\end{enumerate}
which give rise to solutions that propagates to the right side and undergoes dispersion. Since the matrices $\mathbb{M}_j$, $\mathbb{S}_j$ and $\mathbb{M}_{j,W_n}$ are both tri-diagonal symmetric in one dimension,  the consumed memory is low. It is thus possible to solve numerically the Schr{\"o}dinger equation on the entire domain $(0,T) \times (a_0 , b_0)$ with a standard machine. The modulus of solutions at the final time $t = T$ are presented in Figure \ref{Initialsol} for $\mathscr{V}=-x^2$ and $\mathscr{V}=|u|^2$.

We use a cluster consisting of 92 nodes (16 cores/node, Intel Sandy Bridge E5-2670, 32GB/node) to implement the SWR algorithms. We fix one MPI process per subdomain and 16 MPI processes per node. The communications are handled by PETSc and Intel MPI. The linear systems \eqref{AvML} and \eqref{AvMNL} related to the Schr{\"o}dinger equation are solved by the LU direct method using the MKL Pardiso library. The convergence condition for our SWR algorithm is $\parallel g^{k+1} -g^{k} \parallel <10^{-10}$. Two types of initial vectors $g^0$ are considered in this article. One is the zero vector, another is the random vector. According to our tests, the zero initial vector makes the algorithms to converge faster, but obviously it could not include all the frequencies. As mentioned in \cite{Gander2008history}, using the zero initial vector could give wrong conclusions associated with the convergence. Thus, the zero vector is used when one wants to evaluate the computation time, while the random vector is used when comparing the transmission conditions.

%

\subsection{Comparison of classical and new algorithms}
\label{Sec_compare_algorithmes}

We are interested in this part to observe the robustness of the algorithms, to know whether they converge or not for the time independent potential $\mathscr{V}=-x^2$. Similarly, we will observe the computation time and the high scalability of the algorithms. We denote by $T^{\mathrm{ref}}$ the computation time required to solve numerically on a single processor the Schr{\"o}dinger equation on the entire domain and $T^{\mathrm{cls}}$ (resp. $T^{\mathrm{new}}$) the computation time of the classical (resp. new) algorithm for $N$ subdomains. We test the algorithms for $N=2, 10, 100, 500, 1000$ subdomains with the transmission condition $S_0^2$. The reason for using $S_0^2$ for these tests will be explained in Remark \ref{rmk_S02}. The initial vector here is the zero vector.

\begin{figure}[!htbp]
  \centering
    \includegraphics[width=0.46\textwidth]{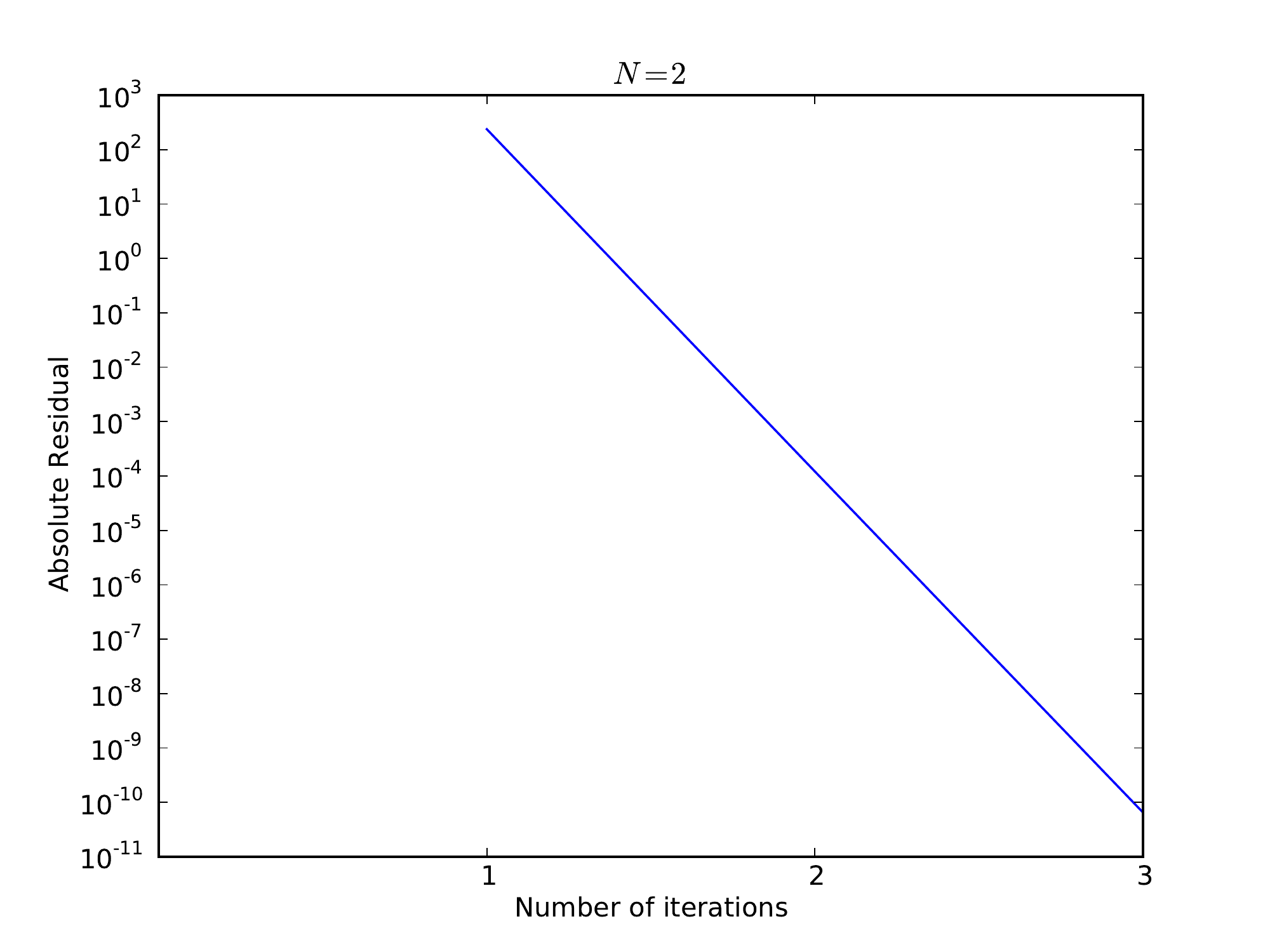}
    \includegraphics[width=0.46\textwidth]{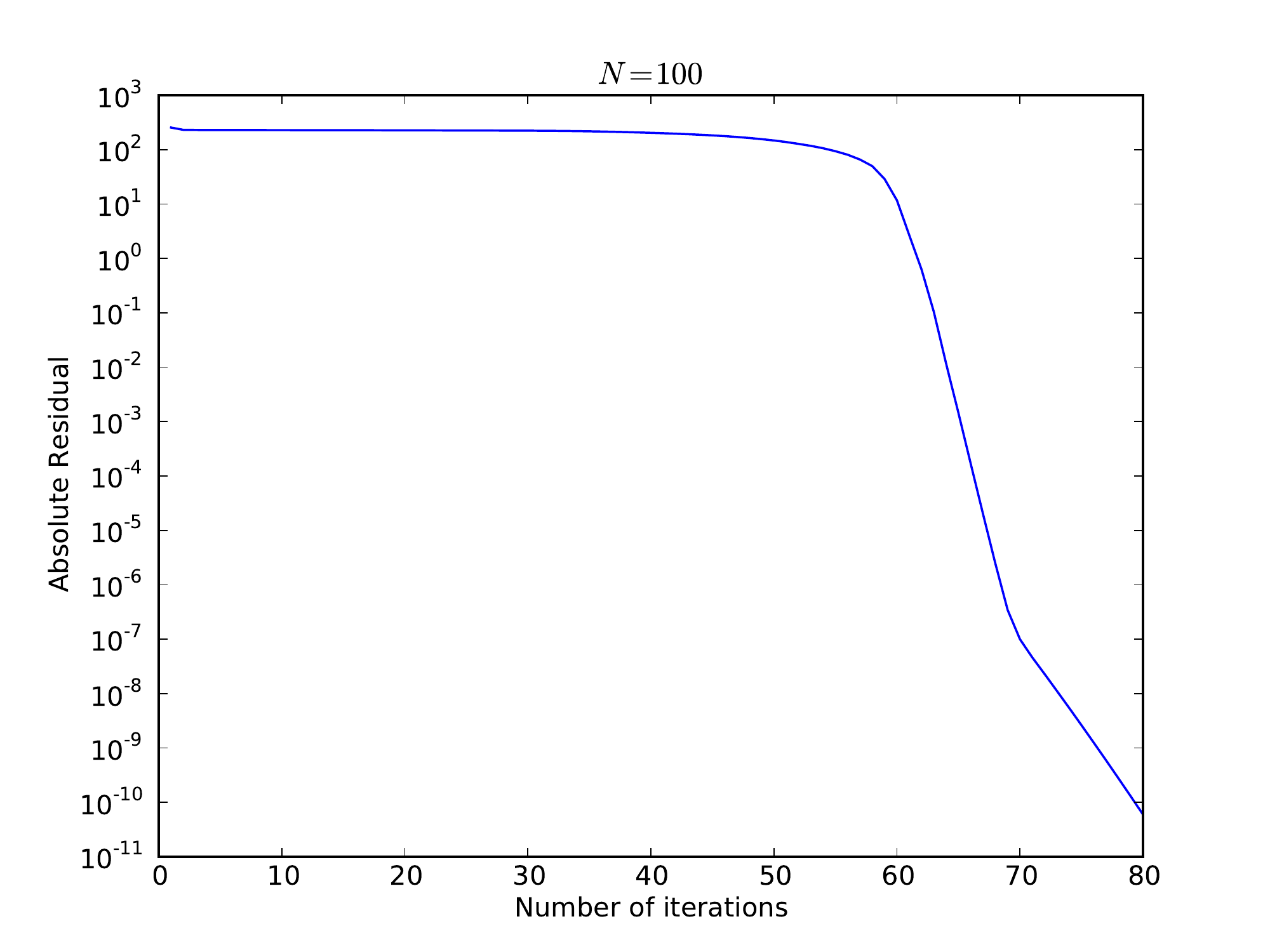}
  \caption{Convergence history, $N=2,100$, $\mathscr{V}=-x^2$, $\Delta t=0.001$, $\Delta x = 10^{-5}$, Fixed point.}
  \label{Fig_S02xx}
\end{figure}

First, the convergence history and the computation time for the algorithms are shown in Figure \ref{Fig_S02xx} and Table \ref{Time_S02_xx} where the fixed point method is used on the interface problem. The algorithms converge for 500 sub domains, but not for 1000 sub domains.

\begin{table}[!htbp]
\footnotesize 
  \centering
   \renewcommand{\arraystretch}{1.5}
    \caption{Computation time in seconds, $\mathscr{V}=-x^2$, $\Delta t=0.001$, $\Delta x = 10^{-5}$, Fixed point.}
  \begin{tabular}{|c|c|c|c|c|}
    \hline
    $N$ & 2 & 10 & 100 & 500 \\
    \hline  
    $T^{\mathrm{ref}}$ & \multicolumn{4}{c|}{403.56} \\
    \hline
    $T^{\mathrm{cls}}$ & 773.07 & 2937.77 & 359.30 & 284.78 \\
    \hline
    $T^{\mathrm{new}}$ & 773.72 & 178.30 & 18.19 & 4.76 \\
    \hline
  \end{tabular}
  \label{Time_S02_xx}
\end{table}

Next, we use the Krylov methods (Gmres or Bicgstab) on the interface problem instead of the fixed point method. Table \ref{Time_S02_xx_gmresbicgstab} present the computation time.
\begin{table}[!htbp] 
\footnotesize 
  \caption{Computation time in seconds, $\mathscr{V}=-x^2$, $\Delta t=0.001$, $\Delta x = 10^{-5}$, Gmres and Bicgstab.}
  \centering
  \renewcommand{\arraystretch}{1.5}
  \begin{tabular}{|c|c|c|c|c|c|c|}
    \hline
    & $N$ & 2 & 10 & 100 & 500 & 1000 \\
    \cline{2-7}  
    & $T^{\mathrm{ref}}$ & \multicolumn{5}{c|}{403.56} \\
    \hline
    \multirow{2}{*}{Gmres}
    & $T^{\mathrm{cls}}$ & 771.82 & 2577.51 & 2249.54 & 907.06 & 739.65 \\
    & $T^{\mathrm{new}}$ & 777.42 & 177.20 & 18.95 & 6.86 & 8.17 \\
    \hline
    \multirow{2}{*}{Bicgstab}
    & $T^{\mathrm{cls}}$ & 774.19 & 2760.11 & 679.72 & 799.09 & 845.65 \\
    & $T^{\mathrm{new}}$ & 774.44 & 177.02 & 18.18 & 6.83 & 7.12 \\
    \hline
  \end{tabular}
  \label{Time_S02_xx_gmresbicgstab}
\end{table}
As we can see, the use of Krylov methods allows to obtain robust scalable SWR algorithms. The algorithms converge for 1000 subdomains and are scalable up to 500 subdomains. Besides their computation times are lower than the ones of the classical algorithm.
%
Roughly speaking, in Table \ref{Time_S02_xx} and Table \ref{Time_S02_xx_gmresbicgstab} we have
\begin{align*}
  & T^{\mathrm{cls}} = T_{\mathrm{sub}} \times N_{\mathrm{iter}} + ...,\\
  & T^{\mathrm{new}} = T_{\mathrm{sub}} \times 4 + T_{Ld} + ...,
\end{align*}
where $T_{\mathrm{sub}}$ is the computation time for solving the equation on one subdomain, $T_{Ld}$ is the computation time for solving the interface problem, ``$...$''  represent the negligible part of computation time such as the construction of matrices for the finite element method. If the number of subdomains $N$ is not so large, then $T_{\mathrm{sub}} \gg T_{Ld}$ and the minimum of $N_{\mathrm{iter}}$ is 3 in all our tests. If the number of subdomains $N$ is large, then $T_{Ld} \sim T_{\mathrm{sub}}$ and $N_{\mathrm{iter}} \gg 4$. It is for this reason that the new algorithm takes less computation time. However, as the number of subdomains increase, $T_{Ld}$ becomes larger. Thus, the new algorithm loses scalability if the number of subdomains is large. 

In conclusion, the new algorithm with Krylov methods is robust and it takes much less computation time than the classical algorithm. 

\subsection{Comparison of classical and preconditioned algorithms}

In this part, we are interested in observing the robustness, the computation time and the scalability of the preconditioned and non-preconditioned (classical) algorithms for time dependent potential $\mathscr{V}=5tx$ and nonlinear potential $\mathscr{V}=|u|^2$. We denote by $N_{\mathrm{pc}}$ the number of iterations required to obtain convergence with the preconditioned algorithm and $T_{\mathrm{pc}}$ the computation time of the preconditioned algorithm. The transmission condition used in this section is $S_0^2$. We use the zero vector as the initial vector $g_0$.

\begin{figure}[!htbp]
  \centering
    \includegraphics[width=0.46\textwidth]{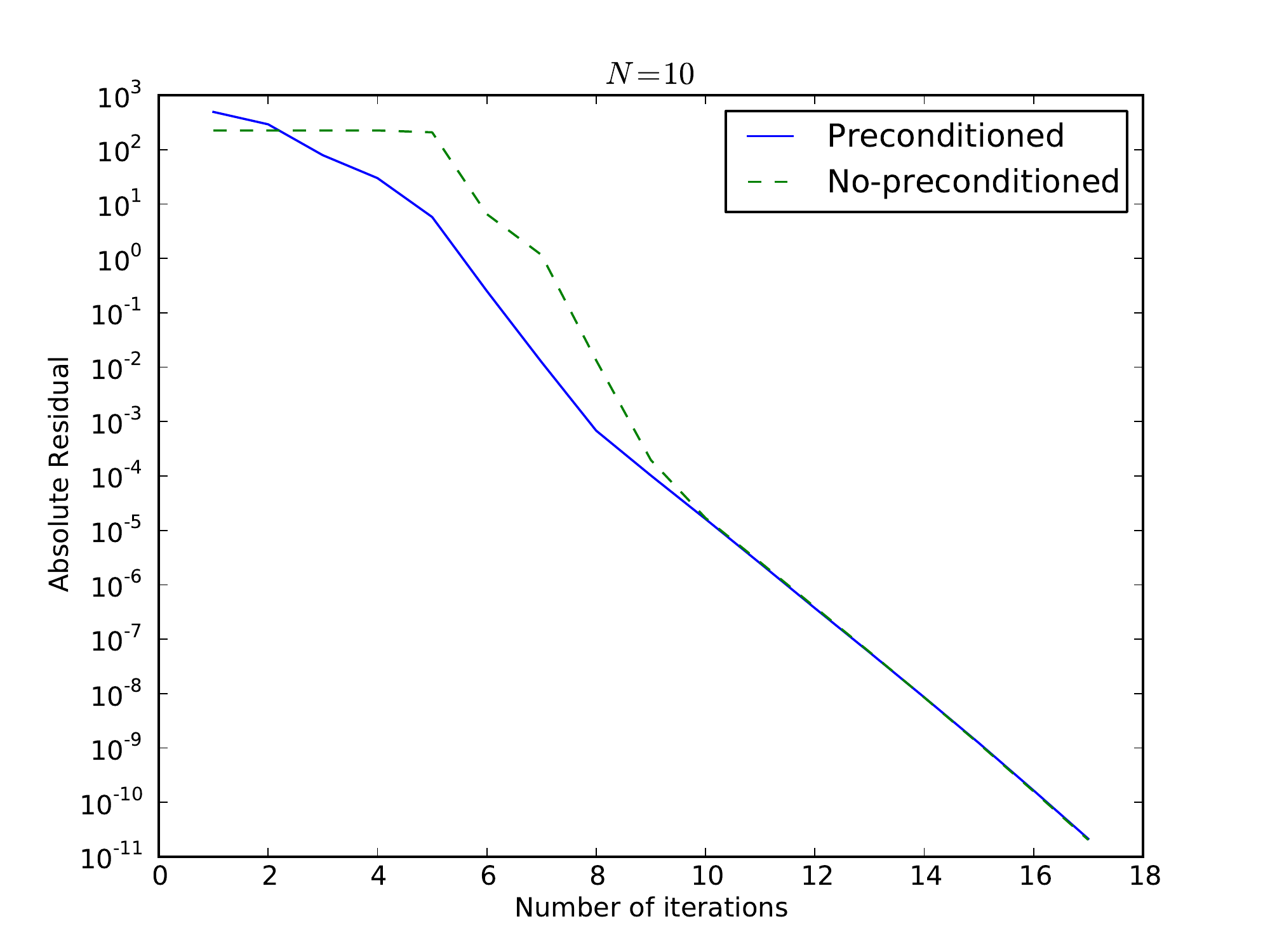}
    \includegraphics[width=0.46\textwidth]{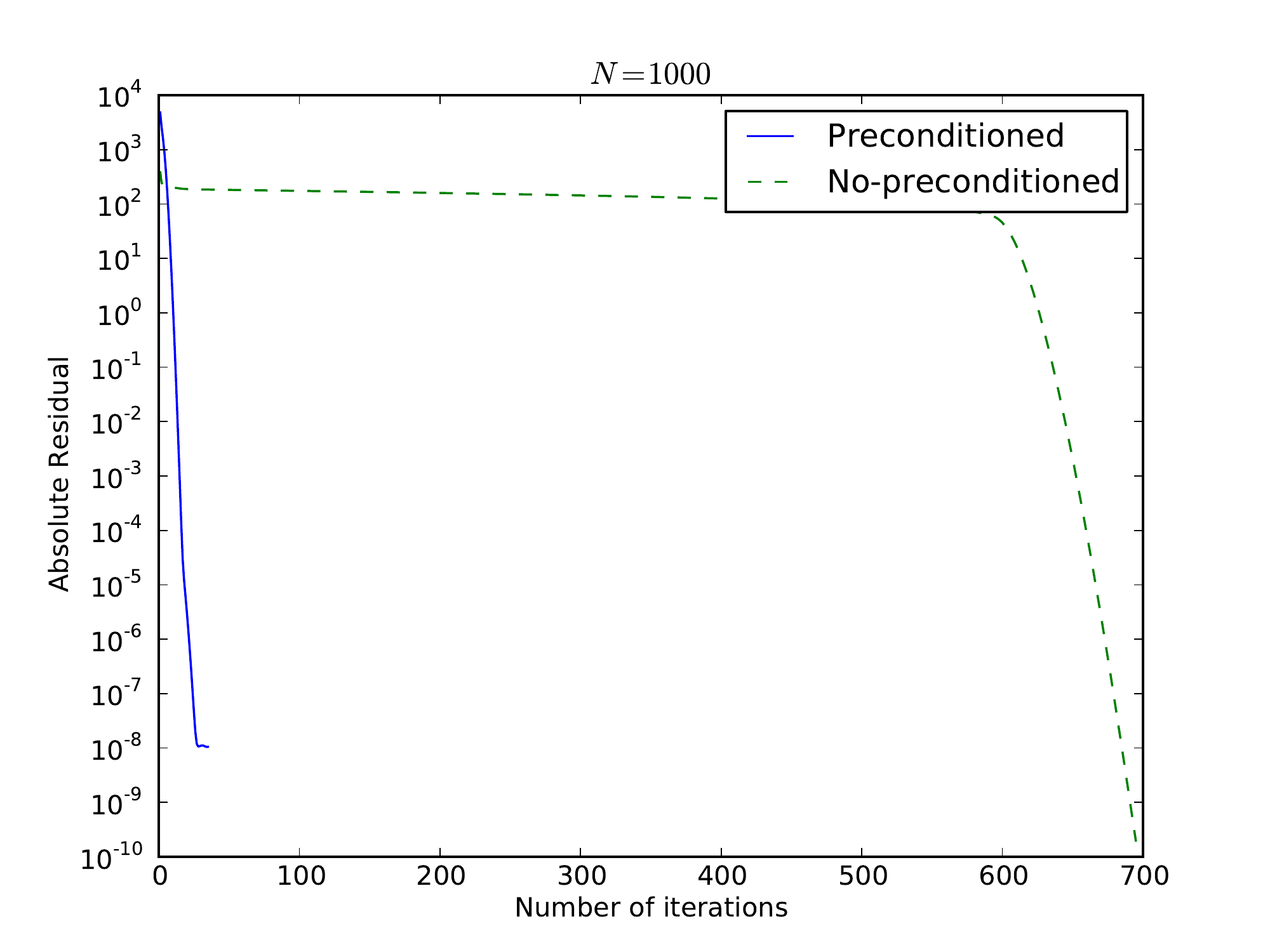}
  \caption{Convergence history, $N=10,1000$,
    $\mathscr{V}=5tx$, $\Delta t = 0.001$, $\Delta x = 10^{-5}$.}
  \label{Fig_hist_Ltx}
\end{figure}

First, we present in Figure \ref{Fig_hist_Ltx} the convergence history for $\mathscr{V}=5tx$. If $N$ is not large, then there is no big difference between the classical algorithm and the preconditioned algorithm. However, if $N$ is large, then as at each iteration, one subdomain communicate only with two adjacent subdomains, we can see that the non-preconditioned algorithm converges very slowly in the first interations. The convergence of the preconditioned algorithm improves greatly since the preconditioner allows communication with remote subdomains.
The number of iterations required for convergence and the computation time are presented in Table \ref{Tab_compa_Ltx_1e-5} for $N=10$, $N=100$, $N=500$ and $N=1000$. 
We can see that the preconditioner allows to decrease significantly both the number of iterations and the computation time. The strong scalability of the classical algorithm is very low. Indeed, the number of iterations required increases with the number of subdomains. The preconditioned algorithm is much more scalable (up to 500 subdomains). However, it loses scalability from $N = 500$ to $N = 1000$. There are two reasons. One is that the number of iterations required for $N = 1000$ is a little bit more than that for $N = 500$. The other one is linked to  the implementation of the preconditioner.
%
Indeed, the time $T_{\mathrm{pc}}$ consists of three major parts: the application of $\mathcal{R}$ to vectors (step 1, denoted by $T_1$), the construction of the preconditioner (denoted by $T_{3c}$) and the application of preconditioner (step 3, denoted by $T_3$). We have thereby
\begin{equation}
  \label{Tpc_123}
  T_{\mathrm{pc}} \approx T_1 + T_{3c} + T_3.
\end{equation}
If $N$ is not very large, $T_1 \sim T_{3c} \gg T_3$. By increasing the number of subdomains, $T_1$ and $T_{3c}$ decreases and $T_3$ increases. Thus, if $N$ is large, $T_3$ is not negligible compared to $T_1$ and $T_{3c}$. However, it is not very convenient to estimate $T_1$ and $T_3$ in our codes because we use the "free-matrix" solvers in the PETSc library. To confirm our explanation, we make tests using a coarser mesh in space ($\Delta t=0.001$, $\Delta x=10^{-4}$). The size of the interface problem \eqref{interfacepbR} is the same, thus $T_3$ should be similar to that of the previous tests ($\Delta t=0.001$, $\Delta x=10^{-5}$). But the size of the problem on a subdomain is ten times smaller. Thus, $T_1$ and $T_{3c}$ are both smaller. The preconditioned algorithm should be less scalable. The results are shown in Table \ref{Tab_compa_Ltx_1e-4}. It can be seen that the computation time $T_{\mathrm{pc}}$ for $N = 1000$ is larger than for $N = 500$ and the preconditioned algorithm is not very scalable from $N=100$ to $N=500$.
Despite this remark, we could conclude from our tests that the preconditioned algorithm reduces a lot the number of iterations and the computing time compared to the classical algorithm.
%
\begin{table}[!htbp]
\footnotesize  
  \centering
    \caption{Number of iterations required and computation time of the classical algorithm and the preconditioned algorithm, $\mathscr{V}=5tx$, $\Delta t = 0.001$, $\Delta x = 10^{-5}$.}
      \renewcommand{\arraystretch}{1.2}
  \begin{tabular}{|c|c|c|c|c|}
    \hline
    $N$ & 10 & 100 & 500 & 1000 \\
    \hline
    $N_{\mathrm{nopc}}$ & 17 & 71 & 349 & 695 \\
    \hline
    $N_{\mathrm{pc}}$ & 17 & 32 & 31 & 35 \\
    \hline
    $T^{\mathrm{ref}}$ & \multicolumn{4}{c|}{6496.3}\\
    \hline
    $T_{\mathrm{nopc}}$ & 10123.1 & 3217.0 & 2466.5 & 2238.0 \\
    \hline
    $T_{\mathrm{pc}}$ & 10128.9 & 1432.7 & 250.0 & 170.7 \\
    \hline
  \end{tabular}
  \label{Tab_compa_Ltx_1e-5}
\end{table}
\begin{table}[!htbp]
\footnotesize  
  \centering
    \caption{Number of iterations required and computation time of the classical algorithm and the preconditioned algorithm, $\mathscr{V}=5tx$, $\Delta t = 0.001$, $\Delta x = 10^{-4}$.}
  \renewcommand{\arraystretch}{1.2}
  \begin{tabular}{|c|c|c|c|c|}
    \hline
    $N$ & 10 & 100 & 500 & 1000 \\
    \hline
    $N_{\mathrm{nopc}}$  & 17 & 71 & 349 & 695 \\
    \hline
    $N_{\mathrm{pc}}$  & 17 & 32 & 26 & 25 \\
    \hline
    $T^{\mathrm{ref}}$ & \multicolumn{4}{c|}{507.5}\\
    \hline
    $T_{\mathrm{nopc}}$  & 681.9 & 223.8 & 210.2 & 191.2 \\
    \hline
    $T_{\mathrm{pc}}$  & 694.3 & 107.6 & 38.4 & 54.5 \\
    \hline
  \end{tabular}
  \label{Tab_compa_Ltx_1e-4}
\end{table}

\begin{figure}[!htbp]
  \centering
    \includegraphics[width=0.46\textwidth]{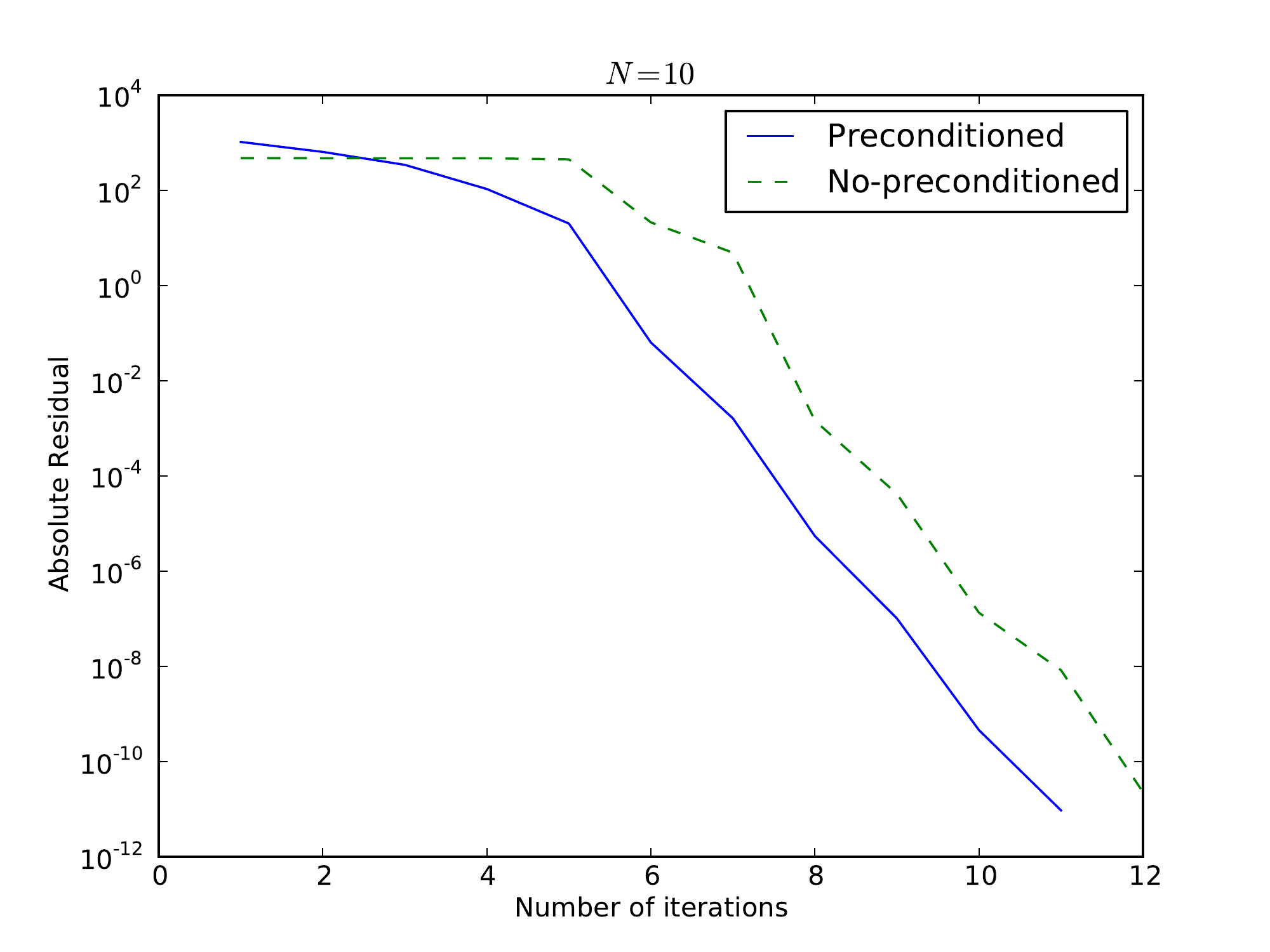}
    \includegraphics[width=0.46\textwidth]{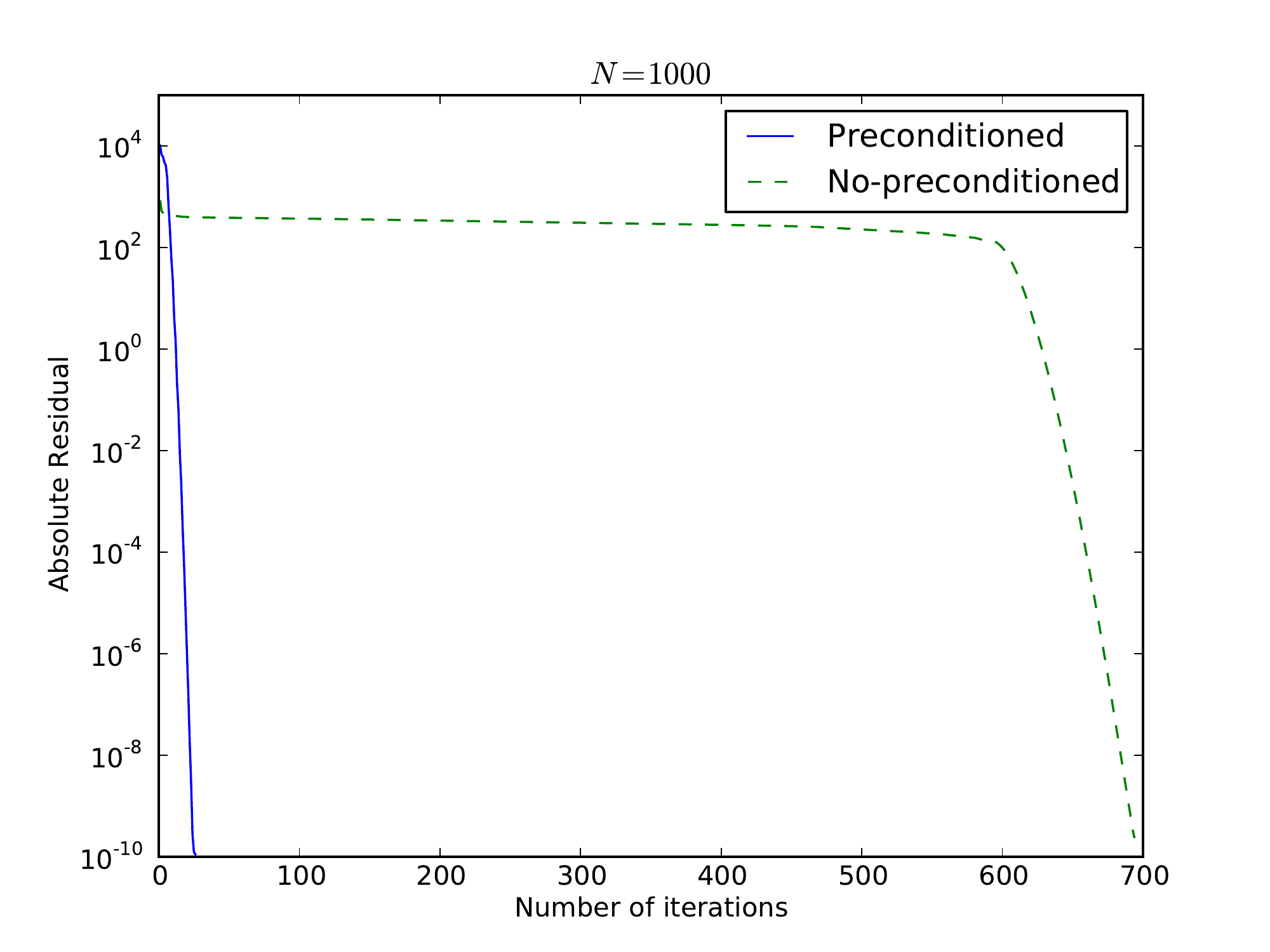}
      \caption{Convergence history, $N=10,1000$,
    $\mathscr{V}=|u|^2$, $\Delta t = 0.001$, $\Delta x = 10^{-5}$.}
    \label{Fig_hist_NL}
\end{figure}

Next, we reproduce the same tests for the nonlinear potential $\mathscr{V}=|u|^2$. The convergence history is presented in Figure \ref{Fig_hist_NL}. We show the number of iterations and the computation time in Table \ref{Tab_compa_NL_1e-5}. The conclusions are quite similar.
%
\begin{table}[!htbp] 
\footnotesize 
  \centering
      \caption{Number of iterations required and computation time of the classical algorithm and the preconditioned algorithm, $\mathscr{V}=|u|^2$, $\Delta t = 0.001$, $\Delta x = 10^{-5}$.}
  \renewcommand{\arraystretch}{1.5}
  \begin{tabular}{|c|c|c|c|c|}
    \hline
    $N$ & 10 & 100 & 500 & 1000 \\
    \hline
    $N_{\mathrm{nopc}}$ & 12 & 71 & 349 & 694 \\
    \hline
    $N_{\mathrm{pc}}$ & 11 & 22 & 25 & 26 \\
    \hline
    $T^{\mathrm{ref}}$ & \multicolumn{4}{c|}{3200.8}\\
    \hline
    $T_{\mathrm{nopc}}$ & 2582.3 & 1332.2 & 1248.0 & 1129.7 \\
    \hline
    $T_{\mathrm{pc}}$ & 2446.7 & 408.2 & 117.6 & 83.8 \\
    \hline
  \end{tabular}
  \label{Tab_compa_NL_1e-5}
\end{table}

\subsection{Comparison of the transmission conditions}
\label{Sec_compare_conditions}

In this part, we compare the transmission conditions which are presented in Section \ref{Sec_SWR} in the framework of the new algorithm for $\mathscr{V}=-x^2$ and the preconditioned algorithm for $\mathscr{V}=|u|^2$. The theoretical optimal parameter $p$ in the transmission condition Robin being not at hand, we seek the best parameter numerically. We use in the subsection the random vector as the initial vector $g_0$ to make sure that all frequencies are present.

\subsubsection{Case of linear potential}

We first consider the linear potential $\mathscr{V} = -x^2$. We compare the number of iterations, the total computation time to perform a complete simulation and the computation time required ($T_{Ld}$) to solve the interface problem in Table \ref{CompareOpsN2V_xx_all} for $N=2$ using the fixed point method, Gmres and Bicgstab methods on the interface problem.
As can be seen, the total computation times are almost identical. The required computation time for solving the interface problem is relatively close to zero compared with the total computation time. Therefore, we are interested rather in the number of iterations. We can make the following observations
%
\begin{enumerate}
\item
  the number of iterations required for the Robin transmission condition
  is greater compared to the other three strategies,
\item
  in each strategy, the number of iterations is not sensitive to order,
\item
  for the Pad{\'e} approximation strategy, the number of iterations decrease as the parameter of Pad{\'e} ($m$) increase.
\end{enumerate}

\begin{table}[!htbp] 
\footnotesize 
  \renewcommand{\arraystretch}{1.5}
  \centering
    \caption{Comparison of transmission conditions for
    $N=2$, $V=-x^2$, $\Delta t = 10^{-3}$, $\Delta x = 10^{-5}$.}
  \label{CompareOpsN2V_xx_all}
  \begin{tabular}{|c|c|c|c|c|c|c|c|c|c|c|}
    \hline
     \multicolumn{2}{|c}{} & \multicolumn{3}{c|}{Fixed point} & \multicolumn{3}{c|}{Gmres} & \multicolumn{3}{c|}{Bicgstab}\\ \hline
    \multicolumn{2}{|c|}{Strategy} & $N_{\mathrm{iter}}$ & $T_{Ld}$ & $T_{\mathrm{total}}$ & $N_{\mathrm{iter}}$ & $T_{Ld}$ & $T_{\mathrm{total}}$ & $N_{\mathrm{iter}}$ & $T_{Ld}$ & $T_{\mathrm{total}}$\\
    \hline
    \multirow{3}{*}{$S_{0}^M$}
    & $S_0^{2}$ & 6 & 0.005 & 775.7 & 5 & 0.002 & 774.2 & 3 & 0.002 & 773.8 \\
    & $S_0^{3}$ & 6 & 0.002 & 774.2 & 5 & 0.002 & 779.6 & 3 & 0.002 & 773.3 \\
    & $S_0^{4}$ & 6 & 0.002 & 769.0 & 5 & 0.002 & 774.2 & 3 & 0.002 & 773.6 \\
    \hline
    \multirow{2}{*}{$S_{1}^M$}
    & $S_1^{2}$ & 6 & 0.002 & 773.4 & 5 & 0.002 & 773.2 & 3 & 0.002 & 773.8\\
    & $S_1^{4}$ & 6 & 0.002 & 773.9 & 5 & 0.002 & 773.6 & 3 & 0.002 & 774.5\\
    \hline
    \multirow{6}{*}{$S_{2}^M$}
    & $S_2^{2,20}$ & 191 & 0.062 & 773.3  & 28 & 0.010 & 774.5 & 16 & 0.011 & 773.1\\
    & $S_2^{2,50}$ & 76  & 0.025 & 773.6  & 27 & 0.010 & 773.3 & 15 & 0.010 & 773.6\\
    & $S_2^{2,100}$& 39  & 0.013 & 776.3  & 23 & 0.008 & 775.2 & 13 & 0.009 & 773.6\\
    & $S_2^{4,20}$ & 181 & 0.059 & 769.9  & 28 & 0.010 & 774.6 & 15 & 0.010 & 773.6\\
    & $S_2^{4,50}$ & 77  & 0.025 & 776.0  & 27 & 0.010 & 773.5 & 15 & 0.010 & 773.3\\
    & $S_2^{4,100}$& 39  & 0.013 & 775.4  & 23 & 0.008 & 773.8 & 13 & 0.009 & 774.8\\
    \hline                                                       
    \multicolumn{2}{|c|}{Robin$^{*}$}  & 1112 & 0.360 & 774.7 & 47 & 0.017 & 776.4 & 27 & 0.018 & 777.4 \\
    \hline
  \end{tabular}
  \begin{tablenotes}
  \item {$*$ the parameters for the transmission condition Robin are $p=44$ (fixed point), $p=5$ (Gmres) and $p=5$ (Bicgstab).}  
  \end{tablenotes}
\end{table}
We make the same tests for $N=500$, the results are shown in Table \ref{CompareOpsN500V_xx_all}. We could see that
\begin{enumerate}
\item
  in each strategy, the number of iterations is not sensitive to order,
  \item
  for the Pad{\'e} approximation strategy, if the parameter $m$ is small, then the algorithm is not robust, 
  \item
    the Krylov methods (Gmres and Bicgstab) could not always reduce the number of iterations.
  \end{enumerate}
  
\begin{table}[!htbp] 
\footnotesize 
  \renewcommand{\arraystretch}{1.5}
  \centering
   \caption{Comparison of transmission conditions for
    $N=500$, $V=-x^2$, $\Delta t = 10^{-3}$, $\Delta x = 10^{-5}$.}
  \begin{tabular}{|c|c|c|c|c|c|c|c|c|c|c|}
    \hline
    \multicolumn{2}{|c|}{} & \multicolumn{3}{c|}{Fixed point} & \multicolumn{3}{c|}{Gmres} & \multicolumn{3}{c|}{Bicgstab}\\ \hline
    \multicolumn{2}{|c|}{Strategy} & $N_{\mathrm{iter}}$ & $T_{Ld}$ & $T_{\mathrm{total}}$ & $N_{\mathrm{iter}}$ & $T_{Ld}$ & $T_{\mathrm{total}}$ & $N_{\mathrm{iter}}$ & $T_{Ld}$ & $T_{\mathrm{total}}$\\
    \hline
    \multirow{3}{*}{$S_0^M$}
    & $S_0^{2}$ & 357 & 0.775 & 4.68 & 1023 & 2.883 & 6.91 & 368 & 1.646 & 5.51 \\
    & $S_0^{3}$ & 337 & 0.734 & 4.62 & 977 & 2.620 & 6.55 & 345 & 1.831 & 5.77 \\
    & $S_0^{4}$ & 337 & 0.733 & 4.65 & 978 & 2.681 & 6.54 & 350 & 1.739 & 5.73 \\
    \hline
    \multirow{2}{*}{$S_1^M$} 
    & $S_1^{2}$ & 341 & 0.745 & 4.62 & 1010 & 2.364 & 6.20 & 353 & 2.102 & 6.00 \\
    & $S_1^{4}$ & 340 & 0.743 & 4.63 & 1023 & 3.454 & 7.19 & 351 & 2.225 & 6.06 \\
    \hline
    \multirow{6}{*}{$S_2^M$}
    & $S_2^{2,20}$ & - & & & 1240 & 3.368 & 7.34 & 440 & 2.626 & 6.64 \\
    & $S_2^{2,50}$ & - & & & 997 & 2.320 & 6.30 & 352 & 2.240 & 6.16 \\
    & $S_2^{2,100}$& 336 & 0.735 & 4.62 & 998 & 3.055 & 7.03 & 333 & 1.603 & 5.62 \\
    & $S_2^{4,20}$ & - & & & 1216 & 3.349 & 7.31 & 464 & 2.044 & 6.05 \\
    & $S_2^{4,50}$ & - & & & 1043 & 3.907 & 7.85 & 336 & 1.756 & 5.63 \\
    & $S_2^{4,100}$& 336 & 0.733 & 4.60 & 1024 & 2.424 & 6.35 & 334 & 1.989 & 5.95 \\
    \hline                                                       
    \multicolumn{2}{|c|}{Robin$^{*}$}  & 1690 & 3.628 & 7.52 & 1060 & 3.000 & 6.80 & 318 & 1.41 & 5.32 \\
    \hline
  \end{tabular}
  \begin{tablenotes}
  \item $*$: the parameters for the transmission condition Robin are $p=45$ (fixed point), $p=19$ (Gmres) and $p=6$ (Bicgstab).
  \item -: the algorithm does not converge before 2000 iterations.
  \end{tablenotes}
  \label{CompareOpsN500V_xx_all}
\end{table}

We could conclude that if the number of subdomains $N$ is not very large, the potential strategy in order 2 with Bicgstab method on the interface problem is a good choice. If $N$ is large, the Bicgstab method also allows most of the algorithms to converge, but it is difficult to have a general conclusion for the transmission conditions in the framework of new algorithm.

\subsubsection{Case of nonlinear potential}

Now we turn to compare the transmission conditions for the nonlinear potential $\mathscr{V}=|u|^2$ in the framework of the preconditioned algorithm. First, we study the influence of the parameter $p$ in the Robin transmission condition. The number of iterations and the computation time are shown in Table \ref{Robin_NL_NiterTime}. It is clear that the convergence is not sensitive to this parameter. 
\begin{table}[!htbp] 
\footnotesize 
  \caption{Influence of parameter $p$ in the transmission conditions Robin, $N=2,10,100$, $\mathscr{V}=|u|^2$, $\Delta t=0.001$, $\Delta x = 10^{-4}$.}
  \renewcommand{\arraystretch}{1.5}
  \centering
  \begin{tabular}{|c|c|c|c|c|c|c|c|}
    \hline
    \multirow{10}{*}{Robin}    
    &  & \multicolumn{2}{c|}{$N=2$} & \multicolumn{2}{c|}{$N=10$} & \multicolumn{2}{c|}{$N=100$} \\
    \hline
    & $p$ & $N_{\mathrm{iter}}$ & $T_{\mathrm{total}}$ & $N_{\mathrm{iter}}$ & $T_{\mathrm{total}}$ &  $N_{\mathrm{iter}}$ & $T_{\mathrm{total}}$ \\
    \hline
    & $5$ & 9 & 1042.6 & 12 & 257.9 & 21 & 55.7 \\
    & $10$ & 8 & 920.3 & 11 & 230.9 & 22 & 50.6 \\
    & $15$ & 8 & 920.3 & 11 & 228.7 & 22 & 46.7 \\
    & $20$ & 8 & 914.5 & 11 & 226.1 & 22 & 43.7 \\
    & $25$ & 8 & 913.0 & 11 & 226.4 & 22 & 43.6 \\
    & $30$ & 8 & 919.2 & 11 & 227.6 & 22 & 43.9 \\
    & $35$ & 8 & 922.1 & 11 & 231.8 & 22 & 44.4 \\
    & $40$ & 8 & 922.8 & 12 & 250.2 & 22 & 45.0 \\
    & $45$ & 8 & 921.7 & 12 & 252.5 & 22 & 46.0 \\
    & $50$ & 8 & 928.3 & 12 & 253.3 & 22 & 46.7 \\
    \hline 
  \end{tabular}
  \label{Robin_NL_NiterTime}
\end{table}
Next we compare the three strategies. The numericals results are presented in Table \ref{Compa_3S_NL}. The transmission conditions $S_0^4$, $S_1^4$ and $S_2^4$ {include the evaluation of $f(u)$}. We don't find a suitable discretization of this term such that the continuity of $\mathbf{v}_j$ at the interfaces ensure the continuity of $\partial_{\mathbf{n}_j} f(u)$. Thus we could not obtain the solution $\mathbf{u}_{j,n}$ that satisfy $\mathbf{u}_{j,n}=R_j \mathbf{u}_{0,n}$. We could see that the number of iterations is not sensitive to the transmission condition and its order. However the computation time for the Pad{\'e} strategy is greater than other strategies. On each subdomain, the non linearity is approximated by a fixed point procedure (see formula \eqref{AvMNL}). This fixed point procedure converges more slowly using the Pad{\'e} strategy than the other strategies. This observation is also found in \cite{Klein2010these}. In conclusion, in the nonlinear case, we also think that the potential strategy of order 2 ($S_0^2$) is a good choice.
\begin{table}[!htbp] 
\footnotesize 
  \caption{Comparison of transmission conditions for $N=2,10,100$, $\mathscr{V}=|u|^2$, $\Delta t=0.001$, $\Delta x = 10^{-4}$.}
  \renewcommand{\arraystretch}{1.5}
  \centering
  \begin{tabular}{|p{1cm}|c|c|c|c|c|c|c|}
    \hline  
    &  & \multicolumn{2}{c|}{$N=2$} & \multicolumn{2}{c|}{$N=10$} & \multicolumn{2}{c|}{$N=100$} \\
    \hline
    & & $N_{\mathrm{iter}}$ & $T_{\mathrm{total}}$ & $N_{\mathrm{iter}}$ & $T_{\mathrm{total}}$ &  $N_{\mathrm{iter}}$ & $T_{\mathrm{total}}$ \\
    \hline
    \multirow{2}{*}{$S_{0}^M$}
    & $S_0^2$ & 8 & 909.5 & 11 & 229.1 & 22 & 40.6 \\
    & $S_0^3$ & 7 & 802.1 & 10 & 205.8 & 22 & 41.6 \\
    \hline 
    \multirow{1}{*}{$S_{1}^M$} 
    & $S_1^2$ & 7 & 802.3 & 10 & 205.6 & 22 & 41.4 \\
    \hline
    \multirow{3}{*}{$S_{2}^M$}
    & $S_2^{2,20}$ & 7 & 1732.5 & 10 & 572.0 & 22 & 128.6 \\
    & $S_2^{2,50}$ & 7 & 4042.9 & 10 & 1342.3 & 23 & 310.3\\
    & $S_2^{2,100}$ & 7 & 7900.5 & 10 & 2640.0 & 22 & 576.0 \\
    \hline
  \end{tabular}
  \label{Compa_3S_NL}
\end{table} 

\begin{rmk}
\label{rmk_S02}
\textnormal{
As we indicated previously, we explain here our choice of transmission condition: the potential strategy of order 2 ($S_0^2$). Indeed, it seems reasonable to consider it since
 \begin{enumerate}
    \item
      the algorithm is robust and the computation time for $S_0^2$  is similar to others transmission conditions,
    \item
      if $N$ is not so large, it is one of the best choice,
    \item
      the implementation of $S_0^2$ is much easier than other transmission conditions.
      \end{enumerate}
    }
\end{rmk}

\subsection{Gpu acceleration}
If the number of subdomain $N$ is not so large, then solving the Schrodinger equation on subdomains takes most of the computation time. We move these computations from Cpu to Gpu. In this subsection, we present the numerical experiments of Gpu acceleration. Two Gpu libraries of NVIDIA are used: CUSPARSE (tri-diagonal solver) and CUBLAS (BLAS operations). We use 8 Gpu Kepler K20, and compare the Cpu and Gpu results for $N=2,4,8$. We use always 1 Gpu/MPI process.
\begin{table}[!htbp] 
\footnotesize
  \centering
    \renewcommand{\arraystretch}{1.5}
    \caption{Cpu and Gpu computation time, Bicgstab, $S_0^2$, $\Delta t=0.001$, $\Delta x=10^{-5}$, $V=-x^2$.\normalsize}
  \begin{tabular}{|c|c|c|c|}
    \hline
    $N$ & 2 & 4 & 8 \\
    \hline
    $T^{\mathrm{Cpu}}$ & 774.4 & 393.0 & 203.2 \\
    \hline
    $T^{\mathrm{Gpu}}$ & 27.90 & 16.13 & 12.54 \\
    \hline
     $T^{\mathrm{Cpu}}/T^{\mathrm{Gpu}}$ & 18 & 24 & 16\\
    \hline
  \end{tabular}
  \label{GpuV_xx1e-5}
\end{table}
Gpu could accelerate a lot the computation as shown in Table \ref{GpuV_xx1e-5}. However the algorithm on Gpu is not scalable. The reason is that the size of problem is not large enough for Gpu. Gpu waste some of its ability. 
We test a larger case only for Gpu: $\Delta t=0.001$, $\Delta x=10^{-6}$. The results are shown in Table \ref{GpuV-xx5e-6}.

\begin{table}[!htbp] 
\footnotesize
\renewcommand{\arraystretch}{1.5}
  \centering
  \caption{Gpu computation time, Bicgstab, $S_0^2$, $\Delta =0.001$, $\Delta x=5\times 10^{-6}$, $V=-x^2$.}
  \begin{tabular}{|c|c|c|c|}
    \hline
    $N$ & 2 & 4 & 8 \\
    \hline
    $T^{\mathrm{Gpu}}$ & 51.95 & 28.21 & 16.30 \\
    \hline
  \end{tabular}
  \label{GpuV-xx5e-6}
\end{table}

Finally, we make the same tests for the nonlinear potential in the framework of the preconditioned algorithm. The results are presented in Table \ref{GpuNL1e-5} and Table \ref{GpuNL5e-6}. The conclusion is similar.

\begin{table}[!htbp]
\footnotesize
\renewcommand{\arraystretch}{1.5}
  \caption{Cpu and Gpu computation time, $\Delta t=0.01$, $\Delta x=10^{-5}$, $V=|u|^2$.}
  \centering
  \begin{tabular}{|c|c|c|c|}
    \hline
    $N$ & 2 & 4 & 8 \\
    \hline
    $T^{\mathrm{Cpu}}$ & 373.6 & 526.7 & 316.0 \\
    \hline
    $T^{\mathrm{Gpu}}$  & 73.9 & 40.1 & 34.0 \\
    \hline
    $T^{\mathrm{Cpu}}/T^{\mathrm{Gpu}}$ & 5 & 13 & 9 \\
    \hline
  \end{tabular}
  \label{GpuNL1e-5}
\end{table}

\begin{table}[!htbp]
\footnotesize
\renewcommand{\arraystretch}{1.5}
  \caption{Gpu computation time, $\Delta t=0.01$, $\Delta x=5\times 10^{-6}$, $V=|u|^2$.}
  \centering
  \begin{tabular}{|c|c|c|c|}
    \hline
    $N$ & 2 & 4 & 8 \\
    \hline
    $T^{\mathrm{Gpu}}$ & 134.3 & 73.7 & 46.0 \\
    \hline
  \end{tabular}
  \label{GpuNL5e-6}
\end{table}

\section{Conclusion}
We proposed in this paper a new algorithm of the SWR method for the one dimensional Schr{\"o}dinger equation with time independent linear potential and a preconditioned algorithm for general potentials. The algorithms for both cases are scalable and could reduce significantly the computation time. Some newly constructed absorbing boundary conditions are used as the transmission condition and compared numerically in the framework of the algorithms proposed by us. We believe that the potential strategy of order 2 is a good choice. Besides, we adapted the codes developed on Cpu to Gpu. According to the experiments, the computation could be accelerated obviously.

\section*{Acknowledgements} We acknowledge Pierre Kestener (Maison de la Simulation Saclay France) for the discussions about the parallel programming, especially for his help about Gpu acceleration. This work was partially supported by the French ANR grant ANR-12-MONU-0007-02 BECASIM (Mod\`eles Num\'eriques call).

\bibliographystyle{elsarticle-num}
\section*{References}
\bibliography{Bib2}

\begin{thebibliography}{10}
\expandafter\ifx\csname url\endcsname\relax
  \def\url#1{\texttt{#1}}\fi
\expandafter\ifx\csname urlprefix\endcsname\relax\def\urlprefix{URL }\fi
\expandafter\ifx\csname href\endcsname\relax
  \def\href#1#2{#2} \def\path#1{#1}\fi

\bibitem{Halpern2010_sch}
L.~Halpern, J.~Szeftel, {Optimized and quasi-optimal Schwarz waveform
  relaxation for the one dimensional Schr\"{o}dinger equation}, Math. Model.
  Methods Appl. Sci. 20~(12) (2010) 2167--2199.

\bibitem{Halpern2006_sch}
L.~Halpern, J.~Szeftel, {Optimized and quasi-optimal Schwarz waveform
  relaxation for the one-dimensional Schr\"{o}dinger equation}, Tech. rep.,
  CNRS (2006).

\bibitem{Caetano2010}
F.~Caetano, M.~J. Gander, L.~Halpern, J.~Szeftel, {Schwarz waveform relaxation
  algorithms for semilinear reaction-diffusion equations}, Networks Heterog.
  Media 5~(3) (2010) 487--505.

\bibitem{Gander2007_diffusion}
M.~J. Gander, L.~Halpern, {Optimized Schwarz Waveform Relaxation Methods for
  Advection Reaction Diffusion Problems}, SIAM J. Numer. Anal. 45~(2) (2007)
  666--697.

\bibitem{Hoang2013}
T.~Hoang, J.~Jaffré, C.~Japhet, M.~Kern, J.~Roberts, {Space-Time Domain
  Decomposition Methods for Diffusion Problems in Mixed Formulations} 51~(6)
  (2013) 3532--3559.

\bibitem{Gander2003_wave}
M.~J. Gander, L.~Halpern, F.~Nataf, {Optimal Schwarz waveform relaxation for
  the one dimensional wave equation}, SIAM J. Numer. Anal. 41~(5) (2003)
  1643--1681.

\bibitem{Halpern2009_wave}
L.~Halpern, J.~Szeftel, {Nonlinear nonoverlapping Schwarz waveform relaxation
  for semilinear wave propagation}, Math. Comput. 78~(266) (2009) 865--889.

\bibitem{Dolean2009max}
V.~Dolean, M.~J. Gander, L.~Gerardo-Giorda, {Optimized Schwarz Methods for
  Maxwell's Equations}, SIAM J. Sci. Comput. 31~(3) (2009) 2193--2213.

\bibitem{Antoine2014}
X.~Antoine, E.~Lorin, A.~D. Bandrauk, {Domain Decomposition Methods and
  High-Order Absorbing Boundary Conditions for the Numerical Simulation of the
  Time Dependent Schr\"{o}dinger Equation with Ionization and Recombination by
  Intense Electric Field}.

\bibitem{Antoine2006}
X.~Antoine, C.~Besse, S.~Descombes, {Artificial boundary conditions for
  one-dimensional cubic nonlinear Schr\"{o}dinger equations}, SIAM J. Numer.
  Anal. 43~(6) (2006) 2272--2293.

\bibitem{Antoine2009a}
X.~Antoine, C.~Besse, P.~Klein, {Absorbing boundary conditions for the
  one-dimensional Schr\"{o}dinger equation with an exterior repulsive
  potential}, J. Comput. Phys. 228~(2) (2009) 312--335.

\bibitem{Antoine2011}
X.~Antoine, C.~Besse, P.~Klein, {Absorbing Boundary Conditions for General
  Nonlinear Schr\"{o}dinger Equations}, SIAM J. Sci. Comput. 33~(2) (2011)
  1008--1033.

\bibitem{Antoine2009}
X.~Antoine, C.~Besse, J.~Szeftel, {Towards accurate artificial boundary
  conditions for nonlinear PDEs through examples}, Cubo, A Math. J. 11~(4)
  (2009) 29--48.

\bibitem{Klein2010these}
P.~Klein, {Construction et analyse de conditions aux limites artificielles pour
  des \'{e}quations de Schr\"{o}dinger avec potentiels et non
  lin\'{e}arit\'{e}s}, Ph.D. thesis, Universit\'{e} Henri Poincar\'{e}, Nancy 1
  (2010).

\bibitem{DUR2000}
A.~Dur\'{a}n, J.~Sanz-Serna, {The numerical integration of relative equilibrium
  solutions. The nonlinear Schrodinger equation}, IMA J. Numer. Anal. 20~(2)
  (2000) 235--261.

\bibitem{Saad2003}
Y.~Saad, {Iterative methods for sparse linear systems}, 2nd Edition, Society
  for Industrial and Applied Mathematics, 2003.

\bibitem{petsc-user-ref}
S.~Balay, M.~F. Adams, J.~Brown, P.~Brune, K.~Buschelman, V.~Eijkhout, W.~D.
  Gropp, D.~Kaushik, M.~G. Knepley, L.~C. McInnes, K.~Rupp, B.~F. Smith,
  H.~Zhang, {PETSc Users Manual}, Tech. Rep. ANL-95/11 - Revision 3.4, Argonne
  National Laboratory (2013).

\bibitem{mpiforum30}
{Message Passing Interface Forum}, {MPI : A Message-Passing Interface Standard
  Version 3.0}, Tech. rep. (2012).

\bibitem{Gander2008history}
M.~J. Gander, {Schwarz methods over the course of time}, Electron. Trans.
  Numer. Anal. 31 (2008) 228--255.

\end{thebibliography}
\end{document}